\documentclass[12pt]{article}
\usepackage{latexsym, amssymb}
\usepackage{dsfont}

\DeclareMathAlphabet\gothic{U}{euf}{m}{n}

\makeatletter
\def\eqnarray{\stepcounter{equation}\let\@currentlabel=\theequation
\global\@eqnswtrue
\tabskip\@centering\let\\=\@eqncr
$$\halign to \displaywidth\bgroup\hfil\global\@eqcnt\z@
  $\displaystyle\tabskip\z@{##}$&\global\@eqcnt\@ne
  \hfil$\displaystyle{{}##{}}$\hfil
  &\global\@eqcnt\tw@ $\displaystyle{##}$\hfil
  \tabskip\@centering&\llap{##}\tabskip\z@\cr}

\def\endeqnarray{\@@eqncr\egroup
      \global\advance\c@equation\m@ne$$\global\@ignoretrue}

\def\@yeqncr{\@ifnextchar [{\@xeqncr}{\@xeqncr[5pt]}}
\makeatother

\begin{document}
\bibliographystyle{tom}

\newtheorem{lemma}{Lemma}[section]
\newtheorem{thm}[lemma]{Theorem}
\newtheorem{cor}[lemma]{Corollary}
\newtheorem{voorb}[lemma]{Example}
\newtheorem{rem}[lemma]{Remark}
\newtheorem{prop}[lemma]{Proposition}
\newtheorem{ddefinition}[lemma]{Definition}
\newtheorem{stat}[lemma]{{\hspace{-5pt}}}

\newenvironment{remarkn}{\begin{rem} \rm}{\end{rem}}
\newenvironment{exam}{\begin{voorb} \rm}{\end{voorb}}
\newenvironment{definition}{\begin{ddefinition} \rm}{\end{ddefinition}}

\newcommand{\gota}{\gothic{a}}
\newcommand{\gotb}{\gothic{b}}
\newcommand{\gotc}{\gothic{c}}
\newcommand{\gote}{\gothic{e}}
\newcommand{\gotf}{\gothic{f}}
\newcommand{\gotg}{\gothic{g}}
\newcommand{\gothh}{\gothic{h}}
\newcommand{\gotk}{\gothic{k}}
\newcommand{\gotm}{\gothic{m}}
\newcommand{\gotn}{\gothic{n}}
\newcommand{\gotp}{\gothic{p}}
\newcommand{\gotq}{\gothic{q}}
\newcommand{\gotr}{\gothic{r}}
\newcommand{\gots}{\gothic{s}}
\newcommand{\gotu}{\gothic{u}}
\newcommand{\gotv}{\gothic{v}}
\newcommand{\gotw}{\gothic{w}}
\newcommand{\gotz}{\gothic{z}}
\newcommand{\gotA}{\gothic{A}}
\newcommand{\gotB}{\gothic{B}}
\newcommand{\gotG}{\gothic{G}}
\newcommand{\gotL}{\gothic{L}}
\newcommand{\gotS}{\gothic{S}}
\newcommand{\gotT}{\gothic{T}}

\newcounter{teller}
\renewcommand{\theteller}{(\alph{teller})}
\newenvironment{tabel}{\begin{list}%
{\rm  (\alph{teller})\hfill}{\usecounter{teller} \leftmargin=1.1cm
\labelwidth=1.1cm \labelsep=0cm \parsep=0cm}
                      }{\end{list}}

\newcounter{tellerr}
\renewcommand{\thetellerr}{(\roman{tellerr})}
\newenvironment{tabeleq}{\begin{list}%
{\rm  (\roman{tellerr})\hfill}{\usecounter{tellerr} \leftmargin=1.1cm
\labelwidth=1.1cm \labelsep=0cm \parsep=0cm}
                         }{\end{list}}

\newcounter{tellerrr}
\renewcommand{\thetellerrr}{(\Roman{tellerrr})}
\newenvironment{tabelR}{\begin{list}%
{\rm  (\Roman{tellerrr})\hfill}{\usecounter{tellerrr} \leftmargin=1.1cm
\labelwidth=1.1cm \labelsep=0cm \parsep=0cm}
                         }{\end{list}}

\newcounter{proofstep}
\newcommand{\nextstep}{\refstepcounter{proofstep}\vertspace \par 
          \noindent{\bf Step \theproofstep} \hspace{5pt}}
\newcommand{\firststep}{\setcounter{proofstep}{0}\nextstep}

\newcommand{\Ni}{\mathds{N}}
\newcommand{\Qi}{\mathds{Q}}
\newcommand{\Ri}{\mathds{R}}
\newcommand{\Ci}{\mathds{C}}
\newcommand{\Ti}{\mathds{T}}
\newcommand{\Zi}{\mathds{Z}}
\newcommand{\Fi}{\mathds{F}}

\newcommand{\proof}{\mbox{\bf Proof} \hspace{5pt}} 
\newcommand{\remark}{\mbox{\bf Remark} \hspace{5pt}}
\newcommand{\vertspace}{\vskip10.0pt plus 4.0pt minus 6.0pt}

\newcommand{\simh}{{\stackrel{{\rm cap}}{\sim}}}
\newcommand{\ad}{{\mathop{\rm ad}}}
\newcommand{\Ad}{{\mathop{\rm Ad}}}
\newcommand{\Aut}{\mathop{\rm Aut}}
\newcommand{\arccot}{\mathop{\rm arccot}}
\newcommand{\capp}{{\mathop{\rm cap}}}
\newcommand{\rcapp}{{\mathop{\rm rcap}}}
\newcommand{\diam}{\mathop{\rm diam}}
\newcommand{\divv}{\mathop{\rm div}}
\newcommand{\codim}{\mathop{\rm codim}}
\newcommand{\RRe}{\mathop{\rm Re}}
\newcommand{\IIm}{\mathop{\rm Im}}
\newcommand{\Tr}{{\mathop{\rm Tr \,}}}
\newcommand{\Vol}{{\mathop{\rm Vol}}}
\newcommand{\card}{{\mathop{\rm card}}}
\newcommand{\supp}{\mathop{\rm supp}}
\newcommand{\sgn}{\mathop{\rm sgn}}
\newcommand{\essinf}{\mathop{\rm ess\,inf}}
\newcommand{\esssup}{\mathop{\rm ess\,sup}}
\newcommand{\Int}{\mathop{\rm Int}}
\newcommand{\lcm}{\mathop{\rm lcm}}
\newcommand{\OPS}{\mathop{\rm OPS}}
\newcommand{\loc}{{\rm loc}}

\newcommand{\mod}{\mathop{\rm mod}}
\newcommand{\spann}{\mathop{\rm span}}
\newcommand{\one}{\mathds{1}}

\hyphenation{groups}
\hyphenation{unitary}

\newcommand{\tfrac}[2]{{\textstyle \frac{#1}{#2}}}
\newcommand{\dfrac}[2]{{\displaystyle \frac{#1}{#2}}}

\newcommand{\ca}{{\cal A}}
\newcommand{\cb}{{\cal B}}
\newcommand{\cc}{{\cal C}}
\newcommand{\cd}{{\cal D}}
\newcommand{\ce}{{\cal E}}
\newcommand{\cf}{{\cal F}}
\newcommand{\ch}{{\cal H}}
\newcommand{\ci}{{\cal I}}
\newcommand{\ck}{{\cal K}}
\newcommand{\cl}{{\cal L}}
\newcommand{\cm}{{\cal M}}
\newcommand{\cn}{{\cal N}}
\newcommand{\co}{{\cal O}}
\newcommand{\cp}{{\cal P}}
\newcommand{\cs}{{\cal S}}
\newcommand{\ct}{{\cal T}}
\newcommand{\cx}{{\cal X}}
\newcommand{\cy}{{\cal Y}}
\newcommand{\cz}{{\cal Z}}

\thispagestyle{empty}

\vspace*{1cm}
\begin{center}
{\Large\bf Analysis of the heat kernel of the  \\[2mm]
Dirichlet-to-Neumann operator} \\[5mm]
\large A.F.M. ter Elst$^1$ and E.M.  Ouhabaz$^2$

\end{center}

\vspace{5mm}

\begin{center}
{\bf Abstract}
\end{center}

\begin{list}{}{\leftmargin=1.8cm \rightmargin=1.8cm \listparindent=10mm 
   \parsep=0pt}
\item
We prove Poisson upper bounds for the kernel $K$ of the semigroup 
generated by the Dirichlet-to-Neumann operator if the underlying domain
is bounded and has a $C^\infty$-boundary.
We also prove Poisson bounds for $K_z$ for all $z$ in the right 
half-plane and  for  all its derivatives. 
\end{list}

\vspace{4cm}
\noindent
February 2013

\vspace{5mm}

\vspace{5mm}

\vspace{15mm}

\noindent
{\bf Home institutions:}    \\[3mm]
\begin{tabular}{@{}cl@{\hspace{10mm}}cl}
1. & Department of Mathematics  & 
  2. & Institut de Math\'ematiques de Bordeaux \\
& University of Auckland   & 
  & Universit\'e Bordeaux 1, UMR 5251,  \\
& Private bag 92019 & 
  &351, Cours de la Lib\'eration  \\
& Auckland 1142 & 
  &  33405 Talence \\
& New Zealand  & 
  & France  \\[8mm]
\end{tabular}

\newpage
\setcounter{page}{1}

\section{Introduction} \label{Spbdton1}

For strongly elliptic operators it is well known that the associated 
semigroup has a kernel which satisfies Gaussian bounds.
On $\Ri^d$ this was proved by Aronson \cite{Aro} and later different proofs were found
to handle operators on domains \cite{Dav2} \cite{Ouh5} \cite{AE1}, 
Laplace--Beltrami operators \cite{Sal3} \cite{Gri5}, subelliptic operators
on Lie groups \cite{VSC} \cite{ER13} \cite{DER4} and references therein.
This subject has attracted attention in the last decades and it is now well understood  that 
Gaussian upper bounds for heat kernels play a fundamental role in problems from harmonic 
analysis such as weak type $(1,1)$ estimates for singular integral operators, 
boundedness of Riesz transforms and  spectral multipliers, $L_p$-analyticity of 
the corresponding  semigroup, 
$L_p$-maximal regularity, $L_p$-independence of the spectrum,\ldots{}.
See Chapter 7 in \cite{Ouh5}  and the monographs mentioned above for an overview on the subject.  

It is our aim in the present paper  to study the heat kernel of the Dirichlet-to-Neumann operator. 
Let $\Omega \subset \Ri^d$ be a bounded connected open set with 
Lipschitz boundary.  
Denote by  $\Gamma = \partial \Omega$  the boundary of $\Omega$, endowed  with the 
$(d-1)$-dimensional Hausdorff measure.
Note that $\Gamma$ is not connected in general.
The Dirichlet-to-Neumann operator $\cn$ is an unbounded operator on $L_2(\Gamma)$ defined as follows. 
Given  $\varphi \in L_2(\Gamma)$,
solve the Dirichlet problem
\begin{eqnarray}
& & \Delta u = 0 \quad \mbox{weakly on } \Omega  \label{eSpbdton1;1}  \\[0pt]
& &  u_{|\Gamma} = \varphi \nonumber
\end{eqnarray}
with $u \in W^{1,2}(\Omega)$.
If $u$ has a weak normal derivative $\frac{\partial u}{\partial \nu}$ in 
$L_2(\Gamma)$, then we say that $\varphi \in D(\cn)$ and 
$\cn \varphi = \frac{\partial u}{\partial \nu}$.
See the beginning of Section~\ref{Spbdton2} for more details on this  definition. 
The Dirichlet-to-Neumann operator, also known as voltage-to-current map,  
arises in the problem  of electrical impedance tomography and in various 
inverse problems (e.g., Calder\'on's problem).  
It is well known that $\cn$ is  positive and self-adjoint, 
so $- \cn$ generates a $C_0$-semigroup $S$ on $L_2(\Gamma)$.
Moreover, $S$ is holomorphic in the right half-plane.
If $\Omega$ has a $C^\infty$-boundary, then $\cn$ is equal to 
$\sqrt{-\Delta_{LB}}$, up to a pseudo-differential operator of order~$0$,
where $\Delta_{LB}$ is the Laplace--Beltrami operator on $\Gamma$
(see Taylor \cite{Tay5} Appendix~C of Chapter~12). 
This implies that $S$ has a smooth kernel $K$.
Since the semigroup generated by $-\Delta_{LB}$ has Gaussian 
kernel bounds, the semigroup generated by $-\sqrt{-\Delta_{LB}}$ satisfies 
Poisson kernel bounds (see, for example, \cite{Yos} page~268).
Therefore one would expect that the kernel of the semigroup $S$ generated by 
$-\cn$ also satisfies Poisson bounds.
It is tempting to use perturbation arguments to achieve this idea but this is highly non-trivial 
because the operators in consideration are not differential operators
(these are pseudo-differential operators).
Nevertheless we shall prove a Poisson upper bound for the heat kernel of $\cn$ and show that 
this is even true for complex time. 
One of the main theorems of this paper reads as follows. 

\begin{thm} \label{tpbdton102}
Suppose $\Omega \subset \Ri^d$ is bounded connected with a $C^\infty$-boundary 
$\Gamma$.
Let $\cn$ be the Dirichlet-to-Neumann map and let $K$ be the kernel of the 
semigroup generated by $-\cn$.
Then there exists a $c > 0$ such that 
\[
|K_z(x,y)|
\leq c \, (\cos \theta)^{-2d(d+1)} \, 
   \frac{(|z| \wedge 1)^{-(d-1)}}{\displaystyle \Big( 1 + \frac{|x-y|}{|z|} \Big)^d }
  \, 
\]
for all $x,y \in \Gamma$ and $z \in \Ci$ with $\RRe z > 0$,
where $\theta = \arg z$.
\end{thm}

We also prove upper  bounds for various derivatives of $K_z$  in Theorem~\ref{tpbdton701}.
As a Corollary of the upper bound with  complex time one obtains  immediately  
that the semigroup generated by 
$-\cn$ on  $L_p(\Gamma)$ is holomorphic on the right half-plane  for all $p \in [1,\infty)$.

For positive time $t$ we prove a more general version of Theorem~\ref{tpbdton102} in which we 
allow  a positive measurable potential.
Let $V \in L_\infty(\Omega)$ and suppose that $V \geq 0$.
Let $\cn_V$ be the Dirichlet-to-Neumann  operator with the condition 
$\Delta u = 0$ in (\ref{eSpbdton1;1}) replaced by $(-\Delta + V) u = 0$ 
weakly on $\Omega$.
Then again $\cn_V$ is a positive self-adjoint operator in $L_2(\Gamma)$ 
(see Section  \ref{Spbdton2}).
We prove the following Poisson bounds for 
the heat kernel of $\cn_V$.

\begin{thm} \label{tpbdton101}
Suppose $\Omega \subset \Ri^d$ is bounded connected with a $C^\infty$-boundary 
$\Gamma$.
Let $V \in L_\infty(\Omega)$ and suppose that 
$V \geq 0$.
Then the semigroup generated by $-\cn_V$ has a kernel $K^V$. 
Moreover, there exists a $c > 0$ such that 
\[
0 \leq K^V_t(x,y)
\leq \frac{c \, (t \wedge 1)^{-(d-1)} \, e^{-\lambda_1 t}}
         {\displaystyle \Big( 1 + \frac{|x-y|}{t} \Big)^d }
\]
for all $x,y \in \Gamma$ and $t > 0$, where $\lambda_1 $ is the first eigenvalue of $\cn_V$.
\end{thm}

The proof of Theorem~\ref{tpbdton101} follows by domination of semigroups.
Indeed, we prove on any Lipschitz domain $\Omega$ that the semigroup
$S^V$ generated by $-\cn_V$ is pointwise dominated by the semigroup $S$.
At first sight, this  is not obvious 
since $\cn_V $ does not seem to be  a perturbation of $\cn$ by some positive potential. 
This domination of semigroups implies the domination of their  corresponding  kernels and  hence 
the Poisson bound for $K^V_t$ follows from that  of $K_t$
for positive  time.
In Section~\ref{Spbdton2} we will prove  positivity and 
domination properties.
Moreover, we prove that the semigroup $S^V$ generated by $-\cn_V$ is 
sub-Markovian and ultracontractive.
This then gives estimates on the $L_p$--$L_q$ norm $\|S^V_t\|_{L_p \to L_q}$ for all $t > 0$ 
and $1 \leq p \leq q \leq \infty$.
These imply the existence of a bounded semigroup kernel for $S^V$ and $S$.
In order to deduce (off-diagonal) Poisson bounds for $S$ we use a
multi-commutator argument of McIntosh and Nahmod \cite{MN}.
If $M_g$ denotes the multiplication operator with a function $g \in C^\infty(\Gamma)$,
then one needs $L_p$--$L_q$  bounds on the commutator $[M_g,S_t]$ and higher order 
commutators $[M_{g},[\ldots,[M_{g},S_t]\ldots]] $. 
Using  Duhamel's  formula these involve commutators like 
$[M_{g},[\ldots,[M_{g},\cn]\ldots]] $,  for which we prove  appropriate $L_p$--$L_q$ bounds 
using   a powerful theorem of Coifman and Meyer \cite{CM2}, 
and Riesz potentials. 
Together with the estimates on $\|S_t\|_{L_p \to L_q}$ for all $t > 0$ 
we then establish Poisson bounds for $K_t$ in Section~\ref{Spbdton4}.
Unfortunately, this proof breaks down if one wants to prove 
Poisson bounds for $K_z$ with $z$ in the right half-plane, since we
do not have appropriate  $L_p$--$L_q$ estimates for  $S_z$.
Nevertheless, using the semigroup $T$ associated to a high enough power
of $\cn$, we will be able, with the Coifman--Meyer commutator bounds, 
Sobolev embedding theorem and spectral theorem, to prove 
bounds on $\|[M_g,T_z]\|_{L_1 \to L_\infty}$ and higher order commutators
in Section~\ref{Spbdton6}. 
By subordination these give bounds for multi-commutators in $S_z$
and then Poisson-type bounds for $K_z$, but with a loss of an~$\varepsilon$.
Luckily, the latter still imply the missing bounds $\|S_z\|_{L_p \to L_q}$
for all $1 \leq p \leq q \leq \infty$.
Then the method in Section~\ref{Spbdton4} gives the bounds of Theorem~\ref{tpbdton102}
for complex $z$.
In Section~\ref{Spbdton7} we deduce Poisson bounds for the 
derivatives of $K_z$.
Finally we discuss holomorphy and $H_\infty$-functional calculus for 
$\cn_V$ and $\cn$ in Section~\ref{Spbdton8}. 
In the appendix we collect definitions and theorems for Sobolev spaces 
on compact manifolds which we need throughout the paper.

Finally, we emphasize that all the methods and heat kernel bounds in this paper are also valid 
if  $\cn$ is the Dirichlet-to-Neumann operator on a compact Riemannian manifold without  boundary. 
In addition all we used is that $\cn$ is a self-adjoint elliptic 
pseudo-differential operator of order~$1$ on a compact Riemannian manifold without boundary. 
Hence one can state  all the results in this setting.

\section{Positivity and domination} \label{Spbdton2}

In this section we define the Dirichlet-to-Neumann operator with a potential.
We then prove that its associated  semigroup on $L_2(\Gamma)$ is sub-Markovian and also prove  domination 
between semigroups associated with Dirichlet-to-Neumann operators with different potentials. 

We assume throughout this section that $\Omega$ is a bounded 
Lipschitz domain of $\Ri^d$. 
(In the rest of this paper we require that $\Omega$ has a $C^\infty$-boundary.)
Let $V \in L_\infty (\Omega,\Ri)$ be a (real-valued) potential.
Define the space $H_V$ of harmonic functions for $-\Delta + V$ by
\[
H_V =\{ u \in W^{1,2}(\Omega):  -\Delta u + V u = 0 \mbox{ weakly on } \Omega \}.
\]
Here and in what follows $-\Delta u + V u = 0$ weakly on $\Omega$ 
means that $u \in W^{1,2}(\Omega)$ and 
\[
\int_\Omega \nabla u . \overline{\nabla \chi} + \int_\Omega V \, u \, \overline \chi  = 0 
\]
for all $\chi \in C_c^\infty(\Omega)$.
Note that we can replace $\chi \in C_c^\infty(\Omega)$ by $\chi  \in W_0^{1,2}(\Omega)$.
Define the continuous sesquilinear form $\gota_V \colon W^{1,2}(\Omega) \times W^{1,2}(\Omega) \to \Ci$ by
\[
\gota_V(u,v) 
= \int_\Omega \nabla u . \overline{\nabla v}
   + \int_\Omega V \, u \, \overline v 
.  \]
It is clear that $H_V$ is a closed subspace of $W^{1,2}(\Omega)$
and 
\begin{equation}
H_V = \{ u \in W^{1,2}(\Omega) : \gota_V(u,v) = 0 \mbox{ for all } v \in \ker \Tr \}
,  
\label{eSpbdton2;1}
\end{equation}
where $\Tr \colon W^{1,2}(\Omega) \to L_2(\Gamma)$ is the trace operator.

Denote by $\Delta_D$ the Laplacian with Dirichlet boundary conditions on 
$\Omega$.
Define the form $\gota_V^D \colon W^{1,2}_0(\Omega) \times W^{1,2}_0(\Omega) \to \Ci$ by
$\gota_V^D = { \gota_V}|_{W^{1,2}_0(\Omega) \times W^{1,2}_0(\Omega)}$.
Then $-\Delta_D + V$ is the operator associated with the form $\gota_V^D$.
If $V \geq 0$, then $0 \notin \sigma(-\Delta_D + V)$.
The space $W^{1,2}(\Omega)$ has the following decomposition.

\begin{lemma} \label{lem2.1} 
Suppose $0 \notin \sigma(-\Delta_D + V)$.
Then
\[
 W^{1,2}(\Omega) = W_0^{1,2}(\Omega) \oplus H_V
.  \]
In particular
\begin{equation}
\Tr(H_V) = \Tr(W^{1,2}(\Omega))
.
\label{eq2.41}
\end{equation}
\end{lemma}
\proof\ 
This result  is already proved in \cite{ArM} Lemma 3.2 when $V$ is constant.
The proof given there  works in our setting but we repeat the arguments for  completeness. 
 
Define $\ca \colon W_0^{1,2}(\Omega) \to W_0^{1,2}(\Omega)'$ by
$\langle \ca u,v \rangle = \gota_V^D(u,v)$.
Since $0 \notin \sigma(-\Delta_D + V)$ it follows from
\cite{ABHN} Proposition~3.10.3 that $\ca$ is invertible.
Let $u \in  W^{1,2}(\Omega)$.
Define $F \in W_0^{1,2}(\Omega)'$ by
\[
F(v) = \int_\Omega \nabla u . \overline{\nabla v}
  + \int_\Omega V \, u \, \overline v
.  \]
Then there exists a unique $u_0 \in W_0^{1,2}(\Omega)$ such that $\ca u_0 = F$.
This means that 
$\langle \ca u_0, \chi \rangle = F(\chi)$ for all $\chi \in W_0^{1,2}(\Omega)$
and hence  
\[
  \int_\Omega \nabla (u - u_0) . \overline{\nabla \chi} 
    + \int_\Omega V \, (u - u_0) \, \overline \chi = 0
.  \]
It follows  that $u - u_0 \in H_V$ and so  $u = u_0 + (u-u_0) \in W_0^{1,2}(\Omega) +  H_V$.
The fact that 
$0 \notin \sigma(-\Delta_D + V)$ implies easily that $W_0^{1,2}(\Omega) \cap H_V = \{0\}$.\hfill$\Box$
 
\vertspace

 A direct  consequence of Lemma \ref{lem2.1} is that the trace $\Tr$ is  
injective as an operator from $H_V$ into $L_2(\Gamma)$.
Indeed, if $u, v \in H_V$ such that 
 $\Tr u = \Tr v $, then $u - v \in H_V \cap W_0^{1,2}(\Omega)$.
Thus $u - v = 0$. 
This is a key ingredient for the next coercivity estimate.

\begin{lemma} \label{lpbdton201}
Suppose $0 \notin \sigma(-\Delta_D + V)$.
Then there are $\mu > 0$ and $\omega \in \Ri$ such that 
\[
\gota_V(u,u) + \omega \, \|\Tr u\|_{L_2(\Gamma)}^2 
\geq \mu \, \|u\|_{W^{1,2}(\Omega)}^2
\]
for all $u \in H_V$.
\end{lemma}
\proof\
Since the embedding of $W^{1,2}(\Omega)$ into $L_2(\Omega)$ is compact,
 it follows that for all $\varepsilon \in (0,1)$ there exists a $c > 0$ such that
\begin{equation} \label{eq2.5}
\int_\Omega | u |^2 
\leq \varepsilon \| u \|_{W^{1,2}(\Omega)}^2 + c \int_\Gamma | \Tr u |^2 
\end{equation}
for all $u \in H_V$.
Therefore, 
\[
\int_\Omega | u |^2 \le \frac{\varepsilon}{1-\varepsilon} \int_\Omega | \nabla u |^2 
+ \frac{c}{1- \varepsilon} \int_\Gamma | \Tr u |^2 
.  \]
Thus
\begin{eqnarray*}
\gota_V(u, u) 
&=& \int_\Omega | \nabla u |^2 + \int_\Omega V \, |u|^2  \\
&\ge& \int_\Omega | \nabla u |^2  
    - \| V \|_\infty  \int_\Omega | u |^2  \\
&\ge&  ( 1 - \frac{\varepsilon}{1-\varepsilon} \| V \|_\infty ) \int_\Omega | \nabla u |^2 
-   \frac{c \, \| V \|_\infty}{1- \varepsilon} \int_\Gamma | \Tr u |^2 .
\end{eqnarray*}
Choosing $\varepsilon = (4(\|V\|_\infty + 1))^{-1}$ one deduces that 
\[
\gota_V(u,u) 
  + \frac{c \, \|V\|_\infty }{1- \varepsilon} \int_\Gamma | \Tr u |^2 
\geq \tfrac{1}{2} \int_\Omega | \nabla u |^2  
.  \]
Hence 
\[
\gota_V(u,u) 
  + \Big( c + \frac{c \, \|V\|_\infty }{1- \varepsilon} \Big) \int_\Gamma | \Tr u |^2 
\geq \tfrac{1}{4} \int_\Omega | \nabla u |^2  + \int_\Omega |u|^2
\geq \tfrac{1}{4} \, \|u\|_{W^{1,2}(\Omega)}^2
\]
by using (\ref{eq2.5}) again.\hfill$\Box$

\vertspace

It follows from (\ref{eSpbdton2;1}) and Lemmas~\ref{lem2.1} and \ref{lpbdton201}
that we can apply \cite{AE2} Corollary~2.2: there exists an
$m$-sectorial operator, which we denote by $\cn_V$, such that for all $\varphi,\psi \in L_2(\Gamma)$
one has $\varphi \in D(\cn_V)$ and $\cn_V \varphi = \psi$ if and only if 
there exists a $u \in W^{1,2}(\Omega)$ such that $\Tr u = \varphi$ and 
\begin{equation}
\int_\Omega \nabla u . \overline{\nabla v}
    + \int_\Omega V \, u  \, \overline v  
= \gota_V(u,v)
= \int_\Gamma  \psi \, \overline{\Tr v} 
\label{elpbdton201;10}
\end{equation}
for all $v \in W^{1,2}(\Omega)$.
Since $\gota_V$ is symmetric, the operator $\cn_V$ is self-adjoint.
Obviously $\cn_V$ is bounded below.
If $\varphi$, $\psi$ and $u$ are as above, then choosing $v \in C_c^\infty(\Omega)$
gives $\Delta u = V \, u \in L_2(\Omega)$ as distribution.
Hence 
\[
\int_\Omega \nabla u . \overline{\nabla v} 
   + \int_\Omega (\Delta u) \, \overline v 
= \int_\Gamma \psi \, \overline{\Tr v}  
\]
for all $ v \in W^{1,2}(\Omega)$ and $\frac{\partial u}{\partial \nu} = \psi$
by the Green formula.
Thus for all $\varphi,\psi \in L_2(\Gamma)$ one has $\varphi \in D(\cn_V)$ 
and $\cn_V \varphi = \psi$ if and only if there exists a $u \in W^{1,2}(\Omega)$
such that $\Tr u = \varphi$, $\Delta u = V \, u$ as distribution and 
$\frac{\partial u}{\partial \nu} = \psi$.

The self-adjoint operator $-\cn_V$ generates a quasi-contraction holomorphic semigroup 
$S^V$ on $L_2(\Gamma)$.  
When $V = 0$ we write for simplicity  $\cn = \cn_0$ and $S = S^0$.

There is another  way to describe the operator $\cn_V$, this time with a form 
with domain in $L_2(\Gamma)$.
Since $\Tr|_{H_V}$ is injective, we can define the form $\gotb_V$ with 
domain $D(\gotb_V) = \Tr(H_V)$ by
\[
\gotb_V(\Tr u, \Tr v) 
= \gota_V(u,v)
\]
for all $u,v \in H_V$.
We equip $D(\gotb_V)$ with the inner product 
$(\Tr u,\Tr v)_{D(\gotb_V)} = (u,v)_{W^{1,2}(\Omega)}$.
Since $H_V$ is closed in $W^{1,2}(\Omega)$ it is clear that 
$D(\gotb_V)$ is a Hilbert space.
It follows from Lemma~\ref{lpbdton201} that the form $\gotb_V$ is 
continuous and elliptic.
Then $\cn_V$ is the operator associated with $\gotb_V$.
Indeed, let $\varphi,\psi \in L_2(\Gamma)$.
Then $\varphi \in D(\cn_V)$ and $\cn_V \varphi = \psi$
if and only if there exists a $u \in W^{1,2}(\Omega)$ such that 
$\varphi = \Tr u$
and (\ref{elpbdton201;10}) is valid for all $v \in W^{1,2}(\Omega)$.
Using (\ref{eSpbdton2;1}) it follows that then $u \in H_V$.
Moreover, if $u \in H_V$, then (\ref{elpbdton201;10}) is valid 
for all $v \in W^{1,2}_0(\Omega)$.
Hence by Lemma~\ref{lem2.1} it is equivalent with the statement 
that there exists a $u \in H_V$ such that $\varphi = \Tr u$ and 
\[
\gotb_V(\varphi,\Tr v) = (\psi, \Tr v)_{L_2(\Gamma)}
\]
for all $v \in H_V$.

In the rest of this section we prove the sub-Markovian property of $S^V$, 
a domination property  and $L_p$--$L_q$ estimates.

\pagebreak[2]

\begin{thm}\label{th2.2} 
\mbox{}
\begin{tabel}
\item \label{th2.2-1} 
If  $-\Delta_D + V \ge 0$ and $ 0 \notin \sigma(-\Delta_D + V)$, 
then the semigroup $S^V$ is positive.
\item \label{th2.2-2} 
If $V \ge 0$ then $S^V$ is sub-Markovian.
\end{tabel}
\end{thm}
\proof\ 
`\ref{th2.2-1}'. 
When $V $ is a constant, the positivity of the semigroup is 
proved in \cite{ArM} Theorem~5.1.
The same proof works here, but we repeat the arguments  for completeness.
By the well known Beurling--Deny criteria (see \cite{Dav2}, Section~1.3
or \cite{Ouh5}, Theorem 2.6), it suffices to prove that 
$\varphi^+ \in D(\gotb_V)$ and $\gotb_V(\varphi^+, \varphi^-) \leq 0$
for all real valued $\varphi \in D(\gotb_V)$.
Let $\varphi \in D(\gotb_V)$ be real valued. 
There exists a $u \in H_V$ such that $\varphi = \Tr u$.
Without loss of generality, $u$ is real valued.
Then $\varphi^+ = \Tr (u^+) \in \Tr(W^{1,2}(\Omega)) = \Tr H_V = D(\gotb_V)$
by (\ref{eq2.41}).
By Lemma \ref{lem2.1} we can write $u^+  = u_0 + u_1 $ and $u^-  = v_0 + v_1 $ 
with $u_0, v_0 \in W_0^{1,2}(\Omega)$ and $u_1, v_1 \in H_V$.
Taking the difference, yields $ u = u^+ - u^- = (u_0 - v_0) + (u_1 - v_1)$.
Since both $u, u_1 -v_1 \in H_V$ it follows that $u_0 = v_0$.
  Therefore with (\ref{eSpbdton2;1}) one deduces that 
\begin{eqnarray*}
\gotb_V (\varphi^+, \varphi^-) 
& = & \gota_V(u_1,v_1)
= \gota_V(u_1, v_0 + v_1)
= \gota_V(u_0 + u_1, v_0 + v_1) - \gota_V(u_0, v_0 + v_1)  \\
& = & \gota_V(u^+,u^-) - \gota_V(u_0, v_0)
= - \gota_V(u_0, v_0)  \\
& = & - \gota_V(u_0, u_0)
= - \int_\Omega (|\nabla u_0|^2 + V \, |u_0|^2)
\leq 0
,
\end{eqnarray*}
since 
\[
\gota_V(u^+,u^-)
= \int_\Omega \nabla (u^+). \nabla(u^-)  + \int_\Omega V \, u^+ \, u^-  = 0
\]
and we used the assumption $-\Delta_D + V \ge 0$ in the last step.
This proves the positivity of the  semigroup $S^V$  on $L_2(\Gamma)$.

`\ref{th2.2-2}'. 
By \cite{Ouh3} or \cite{Ouh5}, Corollary 2.17 it suffices to prove that 
$\one \wedge  \varphi \in D(\gotb_V)$ and $\gotb_V(\one \wedge \varphi, (\varphi - \one)^+) \ge 0$
for all $\varphi \in D(\gotb_V)$ with $\varphi \geq 0$.
Let $\varphi \in D(\gotb_V)$ and suppose $\varphi \geq 0$.
As above, the fact that $\one \wedge  \varphi \in D(\gotb_V)$ follows from (\ref{eq2.41}).
Let $u \in H_V$ be such that $\varphi = \Tr u$.
Without loss of generality, $u$ is real valued.
We decompose $\one \wedge u = u_0 + u_1 \in W_0^{1,2}(\Omega) \oplus H_V$.
Then
\[
 (u-\one)^+ = u - \one \wedge u = (-u_0) + (u-u_1) \in W_0^{1,2}(\Omega) \oplus H_V
.  \]
Using (\ref{eSpbdton2;1}) one deduces that 
\begin{eqnarray*}
\gotb_V(\one\wedge \varphi, (\varphi - \one)^+ ) 
& = & \gota_V(u_1, u - u_1)
= \gota_V(u_0 + u_1, u - u_1)  \\
& = & \gota_V(u_0 + u_1, -u_0 + u - u_1)
   + \gota_V(u_0 + u_1, u_0)  \\
& = & \gota_V(u_0 + u_1, -u_0 + u - u_1)
   + \gota_V(u_0, u_0)  \\
& = & \int_\Omega \nabla (\one \wedge u) .\nabla ((u-\one)^+) + \int_\Omega V (\one \wedge u) \, (u-\one)^+  \\
&& \hspace{60mm} {} +  \int_\Omega  | \nabla u_0 |^2  + \int_\Omega V \, u_0^2   \\
& = & \int_\Omega V \, (u-\one)^+  +  \int_\Omega  | \nabla u_0 |^2  + \int_\Omega V \, u_0^2 
\geq 0
\end{eqnarray*}
as required.\hfill$\Box$

\vertspace

Note that the second part of the previous result can also be deduced from the next theorem in which we prove the domination property. 

\begin{thm}\label{th2.3}
Let $V_1, V_2 \in L_\infty(\Omega, \Ri)$ be such that 
$V_1 \le V_2$, $-\Delta_D + V_1 \ge 0$ and $0 \notin \sigma(-\Delta_D + V_1)$.
Then 
\[
0 \le S_t^{V_2} \varphi \le S_t^{V_1} \varphi 
\]
pointwise for all $t > 0$ and  $0 \le  \varphi \in L_2(\Gamma)$.
In particular, if  $0 \le V \in L_\infty(\Omega)$,  then
\[
0 \le S_t^V \varphi \le S_t \varphi
\]
for all $t > 0$ and $0 \le  \varphi \in L_2(\Gamma)$.
\end{thm}
\proof\   
Using criteria for domination of semigroups (see \cite{Ouh3} or \cite{Ouh5}, Theorem 2.24) 
it suffices to prove that 
\begin{equation}\label{eq2.8}
\gotb_{V_2}(\varphi, \psi) \ge \gotb_{V_1}(\varphi, \psi)
\end{equation}
for all $0 \le \varphi, \psi \in D(\gotb_{V_1})$.
Note that 
\[
D(\gotb_{V_1}) = \Tr (W^{1,2}(\Omega)) = D(\gotb_{V_2})
\]
and the ideal property in \cite{Ouh3} or \cite{Ouh5} is satisfied since both 
semigroups $S^{V_1}$ and $S^{V_2}$ are positive by Theorem \ref{th2.2}
(see Proposition 2.20 in \cite{Ouh5}).

Let $0 \le \varphi, \psi \in D(\gotb_{V_1})$.
There are real valued $u_1,v_1 \in H_{V_1}$ and $u_2,v_2 \in H_{V_2}$
such that $\Tr u_1 = \Tr u_2 = \varphi$ and $\Tr v_1 = \Tr v_2 = \psi$.
Since $u_2 - u_1 \in W^{1,2}_0(\Omega)$ and $v_2 \in H_{V_2}$ one has 
\begin{eqnarray*}
\gotb_{V_2}(\varphi, \psi) 
&=& \gota_{V_2}(u_2, v_2)
= \gota_{V_2}(u_1, v_2)  \\
& = & \gota_{V_1}(u_1, v_2) + \int_\Omega (V_2 - V_1) \, u_1 \, v_2  \\
& = & \gota_{V_1}(u_1, v_1) + \int_\Omega (V_2 - V_1) \, u_1 \, v_2
= \gotb_{V_1}(\varphi, \psi)  + \int_\Omega (V_2 - V_1) \, u_1 \, v_2
.
\end{eqnarray*}
By the lemma below, we show that $u_1 \geq 0$ and $v_2 \geq 0$.
Hence $\int_{\Omega} (V_2 - V_1) \, u_1 \, v_2  \ge 0$ and (\ref{eq2.8}) follows.\hfill$\Box$

\vertspace

We have the following maximum principle. 

\begin{lemma}\label{lem2.4} 
Suppose that $V \in L_\infty(\Omega, \Ri)$ with $-\Delta_D + V \ge 0$ and  
$0 \notin \sigma(-\Delta_D + V)$.
Let $0 \le  \varphi \in D(\gotb_V)$ and let $u \in H_V$ be real valued such that 
$\Tr u = \varphi$.
Then $u \geq 0$ on $\Omega$.
\end{lemma}
\proof\ 
By definition of $u \in H_V$ one has
\[
\int_\Omega \nabla u . \nabla \chi  + \int_\Omega V \, u \, \chi = 0
\]
for all $\chi \in W_0^{1,2}(\Omega)$.
 Note that $u^- = 0$ on $\Gamma$ since $u = \varphi \ge 0 $ on $\Gamma$.
 Hence $u^- \in W_0^{1,2}(\Omega)$ by \cite{Alt} Lemma~A.6.10
and we can choose $\chi = u^-$.
We obtain
\[
\int_\Omega \nabla u . \nabla (u^-) + \int_\Omega V \, u \, u^- = 0
.  \]
Because $\int_\Omega \nabla(u^+). \nabla(u^-)  = 0$ we  arrive at 
\[
\int_\Omega | \nabla (u^-) |^2  + \int_\Omega V \, |u^-|^2  = 0
.  \]
Since  $-\Delta_D + V \ge 0$ and 
$ 0 \notin \sigma(-\Delta_D + V)$ we conclude that $u^- = 0$.\hfill$\Box$

\vertspace

Now we prove $L_p$--$L_q$ estimates for the semigroup $S^V$.
Note that $\lambda_1 \geq 0$ in the next theorem.

\begin{thm}\label{th2.5}
Suppose that $d \ge 2$, let $0 \le V \in L_\infty(\Omega)$ and let
$\lambda_1 \in \sigma(\cn_V)$ be the first eigenvalue of $\cn_V$.
Then for all $1 \le p \le q \le \infty$ and $t > 0$ the operator 
$S_t^V$ is bounded from $L_p(\Omega)$ into $L_q(\Omega)$.
Moreover, there exists a $C > 0$ such that 
\[
\|S_t^V \|_{p \to q} 
\le C \, (t \wedge 1)^{-(d-1)(\frac{1}{p} - \frac{1}{q})} \, e^{-\lambda_1 t} 
\]
for all $t > 0$ and $p,q \in [1,\infty]$ with $p \leq q$.
\end{thm}
\proof\ Suppose first that $d \ge 3$.
By Theorem 2.4.2 in \cite{Nec2}, the trace $\Tr$ is a bounded operator from $D(\gotb_V)$ into
$L_s(\Gamma)$, where $ s = \frac{2(d-1)}{d-2}$.
This implies that there exists a $C \geq 1$ such that 
\[
 \| S_t^V \varphi \|_s^2 \leq C (\gotb_V(S_t^V \varphi, S_t^V \varphi) + \| S_t^V \varphi \|_2^2)
\]
for all $t > 0$ and $\varphi \in L_2(\Gamma)$.
Therefore, $S_t^V$ maps $L_2(\Gamma)$ into $L_s(\Gamma)$ with 
\[
\|S_t^V \|_{2 \to s} \leq C \, t^{-1/2} \, e^t
.  \]
Since the semigroup $S_t^V$ is sub-Markovian by Theorem \ref{th2.2}, 
the last estimate extrapolates  and provides the $L_1$--$L_\infty$ estimate 
\[
\|S_t^V \|_{1 \to\infty} 
\le C' \, t^{-(d-1)} \, e^t
\]
for a suitable $C' > 0$, uniformly for all $t > 0$,
see \cite{Cou2} or \cite{Ouh5}, Lemma~6.1.
By \cite{Ouh5}, Lemma~6.5, the last  estimate improves to
\[
 \|S_t^V \|_{1 \to \infty} 
\le C'' \, t^{-(d-1)} \, e^{-\lambda_1 t} \, (1 + t)^{d-1}.
\]
The conclusion of the theorem follows by interpolation.

If $d = 2$, we apply  the same arguments and use Theorem 2.4.6 in \cite{Nec2}.\hfill$\Box$

\section{Smoothing properties for commutators} \label{Spbdton3}

Let $(M,g)$ be a compact Riemannian manifold (without boundary) of dimension $m$.
For general definitions and theorems on compact Riemannian manifolds
we refer to the appendix.
We emphasize that we do not assume that $M$ is connected.
Then $M$ has a finite number of connected components, say 
$M_1,\ldots,M_N$, with $M_i \neq M_j$ if $i \neq j$.
For all $i \in \{ 1,\ldots,N \} $ the component $M_i$ is a 
compact connected Riemannian manifold.
Therefore it has a natural Riemannian distance, denoted by $d_{M_i}$.
We denote by $\diam M_i$ its diameter.
Set $D = 1 + \sum_{i=1}^N \diam M_i$.
We wish to define a distance on the full manifold.
For all $i \in \{ 1,\ldots,N \} $ fix once and for all an element
$x_i \in M_i$.
Let 
\begin{equation}
W = \{ g \in C^\infty(M,\Ri) : 
          \max_{i,j \in \{ 1,\ldots,N \} } |g(x_i) - g(x_j)| 
          + D \, \|\nabla g\|_\infty \leq D 
    \}
.  
\label{eSpbdton3;1}
\end{equation}
If $x,y \in M$ and $g \in W$, then there are $i,j \in \{ 1,\ldots,N \} $ such 
that $x \in M_i$ and $y \in M_j$.
Note that $\|\nabla(g|_{M_i})\|_{L_\infty(M_i)} \leq 1$.
Therefore $|g(x) - g(x_i)| \leq d_{M_i}(x,x_i) \leq \diam M_i$.
Similarly, $|g(y) - g(x_j)| \leq \diam M_j$.
Moreover, $|g(x_i) - g(x_j)| \leq D$.
Hence $|g(x) - g(y)| \leq 3D$.
Since this is for all $g \in W$, we can define the function
$\rho_M \colon M \times M \to [0,\infty)$ by 
\begin{equation}
\rho_M(x,y)
= \sup \{ |g(x) - g(y)| : g \in W \} 
.  
\label{eSpbdton3;2}
\end{equation}
We collect some properties of $\rho_M$.

\begin{lemma} \label{lpbdton300.5}
\mbox{}
\begin{tabel}
\item \label{lpbdton300.5-1}
The function $\rho_M$ is a metric on $M$, bounded by $3D$.
\item \label{lpbdton300.5-2}
If $i \in \{ 1,\ldots,N \} $ then $\rho_M|_{M_i \times M_i} = d_{M_i}$.
\item \label{lpbdton300.5-3}
If $i,j \in \{ 1,\ldots,N \} $, $x \in M_i$, $y \in M_j$ and $i \neq j$,
then $\rho_M(x,y) \geq 1$.
\item \label{lpbdton300.5-4}
Suppose $k \in \Ni$ and $M$ is embedded in $\Ri^k$.
Then there exists a $c > 0$ such that 
\[
c^{-1} \, |x-y| \leq \rho_M(x,y) \leq c \, |x-y|
\]
for all $x,y \in M$.
\end{tabel}
\end{lemma}
\proof\
Clearly $\rho_M$ satisfies the triangle inequality and is symmetric.
If $i \in \{ 1,\ldots,N \} $ and $x,y \in M_i$, then 
$\rho_M(x,y) \leq d_{M_i}(x,y)$.
Conversely, if $\tilde g \in C^\infty(M_i,\Ri)$ and $\|\nabla \tilde g\|_{L_\infty(M_i)} \leq 1$
then one can define $g \in C^\infty(M,\Ri)$ by $g(z) = \tilde g(z)$ if $z \in M_i$ 
and $g(z) = \tilde g(x_i)$ if $z \not\in M_i$.
Then $g \in W$ and $|\tilde g(x) - \tilde g(y)| = |g(x) - g(y)| \leq \rho_M(x,y)$.
Hence $d_{M_i}(x,y) \leq \rho_M(x,y)$.
Therefore $\rho_M|_{M_i \times M_i} = d_{M_i}$.
Finally, let $i,j \in \{ 1,\ldots,N \} $, $x \in M_i$ and $y \in M_j$
with $i \neq j$.
It is easy to see that $\rho_M(x_i,x_j) \geq D$.
Hence $\rho_M(x,y) \geq D - \rho_M(x,x_i) - \rho_M(y,x_j) \geq 1$.
The last statement follows from Lemma~\ref{lpbdtonA06} and the fact that 
the compact components $M_i$ are disjoint.\hfill$\Box$

\vertspace

Although we do not need the following definition until Section~\ref{Spbdton6},
it is convenient to state it now.
Let $k \in \Ni$.
Define 
\begin{equation}
W_k = \{ g \in C^\infty(M,\Ri) : 
          \max_{ i,j \in \{ 1,\ldots,N \} } |g(x_i) - g(x_j)|  
         + D \, \max_{ \ell \in \{ 1,\ldots,k \} } \|\nabla^\ell g\|_\infty   
          \leq D 
      \}
.  
\label{eSpbdton3;3}
\end{equation}
Clearly $W_1 \supset W_2 \supset \ldots$.
Define 
$\rho^{(k)}_M \colon M \times M \to [0,\infty)$ by 
\[
\rho^{(k)}_M(x,y)
= \sup \{ |g(x) - g(y)| : g \in W_k \} 
.  
\]
Then $\rho^{(1)}_M(x,y) \geq \rho^{(2)}_M(x,y) \geq \ldots$.

\begin{lemma} \label{lpbdton300.6}
Let $k \in \Ni$.
The function $\rho^{(k)}_M$ is a metric on $M$ and it is equivalent to
$\rho_M$.
\end{lemma}
\proof\
Note that for all $i \in \{ 1,\ldots,N \} $ the map
\[
(x,y) \mapsto \sup \{ |g(x) - g(y)| : g \in C^\infty(M_i) \mbox{ and } 
      \|\nabla^\ell g\|_\infty \mbox{ for all } \ell \in \{ 1,\ldots,k \} \}
\]
is a metric on $M_i$ which is equivalent to $d_{M_i}$.
(See Lemma~\ref{lpbdtonA05}.)
Then the first part of the lemma follows as in the proof of Lemma~\ref{lpbdton300.5}.
Moreover, the second part follows from this equivalence.\hfill$\Box$

\vertspace

In the proofs we need various estimates on commutators of pseudo-differential operators
with $C^\infty(M)$-functions.
On $\Ri^m$ these read as follows.
We denote by $\cs(\Ri^m)$ the Schwartz space.

\begin{thm} \label{tpbdton301}
Let $k \in \Ni$ and $T \in \OPS^k(\Ri^m)$.
Let $n \in \{ k,\ldots,k+m \} $.
\begin{tabel}
\item \label{tpbdton301-1}
If $n = k$ then for all $p \in (1,\infty)$ there exists a $c > 0$ 
such that 
\[
\|[M_{g_1},[\ldots,[M_{g_n},T]\ldots]] u\|_p 
\leq c \, \|\nabla g_1\|_\infty \ldots \|\nabla g_n\|_\infty \, \|u\|_p
\]
for all $g_1,\ldots,g_n \in \cs(\Ri^m)$ and $u \in C_c^\infty(\Ri^m)$.
\item \label{tpbdton301-2}
If $n \in \{ k+1,\ldots,k+m-1 \} $ 
then for all $p \in (1,\frac{m}{n-k})$ there exists a $c > 0$ such that 
\[
\|[M_{g_1},[\ldots,[M_{g_n},T]\ldots]] u\|_q
\leq c \, \|\nabla g_1\|_\infty \ldots \|\nabla g_n\|_\infty \, \|u\|_p
\]
for all $g_1,\ldots,g_n \in \cs(\Ri^m)$ and $u \in C_c^\infty(\Ri^m)$,
where $\frac{1}{p} - \frac{1}{q} = \frac{n-k}{m}$.
\item \label{tpbdton301-3}
If $n = k+m$ 
then there exists a $c > 0$ such that 
\[
\|[M_{g_1},[\ldots,[M_{g_n},T]\ldots]] u\|_\infty
\leq c \, \|\nabla g_1\|_\infty \ldots \|\nabla g_n\|_\infty \, \|u\|_1
\]
for all $g_1,\ldots,g_n \in \cs(\Ri^m)$ and $u \in C_c^\infty(\Ri^m)$.
\end{tabel}
\end{thm}
\proof\
Statement~\ref{tpbdton301-1} follows from \cite{CM2} Th\'eor\`eme~2.

Next suppose that $n \in \{ k+1,\ldots,k+m \} $.
Let $K$ be the (distributional) kernel of $T$.
Since $T \in OPS^k(\Ri^m)$, there exists a $c > 0$ such that 
$|K(x,y)| \leq c \, |x-y|^{-m-k}$ for all $x,y \in \Ri^m$ with $x \neq y$.
(See \cite{Ste3} Proposition~VI.4.1.)
Let $g_1,\ldots,g_n \in \cs(\Ri^m)$.
Let $\widetilde K$ denote the kernel of $[M_{g_1},[\ldots,[M_{g_n},T]\ldots]]$.
Then $\widetilde K(x,y) = K(x,y) \prod_{j=1}^n (g_j(x) - g_j(y))$ for all $x \neq y$.
Hence 
\[
|\widetilde K(x,y)| \leq \frac{c \, \|\nabla g_1\|_\infty \ldots \|\nabla g_n\|_\infty}{|x-y|^{m-(n-k)}}
\]
for all $x,y \in \Ri^m$ with $x \neq y$.

If $n \in \{ k+1,\ldots,k+m-1 \} $ then $|\widetilde K|$ is a 
Riesz potential and the boundedness of the multi-commutator from $L_p$ into $L_q$
follows from \cite{Ste1} Theorem~V.1.

Finally, if $n = k+m$ then $\widetilde K$ is bounded. 
Therefore the multi-commutator is bounded from $L_1$ into $L_\infty$.\hfill$\Box$

\vertspace

The theorem transfers to compact Riemannian manifolds.
We emphasize that the manifold does not have to be connected in the 
next proposition.

\begin{prop} \label{ppbdton302}
Suppose $M$ is compact.
Let $k \in \Ni$ and $T \in \OPS^k(M)$.
Let $n \in \{ k,\ldots,k+m \} $.
\begin{tabel}
\item \label{ppbdton302-1}
If $n = k$ then for all $p \in (1,\infty)$ there exists a $c > 0$ 
such that 
\[
\|[M_{g_1},[\ldots,[M_{g_n},T]\ldots]] u\|_p 
\leq c \, \|u\|_p
\]
for all $u \in C^\infty(M)$ and $g_1,\ldots,g_n \in W$.
\item \label{ppbdton302-2}
If $n \in \{ k+1,\ldots,k+m-1 \} $ 
then for all $p \in (1,\ldots,\frac{m}{n-k})$ there exists a $c > 0$ such that 
\[
\|[M_{g_1},[\ldots,[M_{g_n},T]\ldots]] u\|_q
\leq c \, \|u\|_p
\]
for all $u \in C^\infty(M)$ and $g_1,\ldots,g_n \in W$,
where $\frac{1}{p} - \frac{1}{q} = \frac{n-k}{m}$.
\item \label{ppbdton302-3}
If $n = k+m$ 
then there exists a $c > 0$ such that 
\[
\|[M_{g_1},[\ldots,[M_{g_n},T]\ldots]] u\|_\infty
\leq c \, \|u\|_1
\]
for all $u \in C^\infty(M)$ and $g_1,\ldots,g_n \in W$.
\end{tabel}
\end{prop}
\proof\
Since $M$ is compact there are $L \in \Ni$ and for all 
$\ell \in \{ 1,\ldots,L \} $ there exist an
open $U_\ell \subset M$, a $C^\infty$-diffeomorphism
$\varphi_\ell \colon U_\ell \to B(0,1)$  and 
$\chi_\ell, \widetilde \chi_\ell \in C_c^\infty(U_\ell)$ 
such that $\sum_{\ell=1}^L \chi_\ell = \one$
and $\widetilde \chi_\ell(x) = 1$ for all $x \in \supp \chi_\ell$.
Without loss of generality we may assume that there exists a $c_0 > 0$
such that 
\begin{eqnarray*}
\|\nabla(g \circ \varphi_\ell^{-1})\|_{L_\infty(B(0,1))}
& \leq & c_0 \, \|\nabla g\|_{L_\infty(M)} ,  \\
\|u\|_{L_q(U_\ell)}
& \leq & c_0 \, \|u \circ \varphi_\ell^{-1}\|_{L_q(B(0,1))} \mbox{ and }  \\
\|v \circ \varphi_\ell^{-1}\|_{L_p(B(0,1))}
& \leq & c_0 \, \|v\|_{L_p(U_\ell)} 
\end{eqnarray*}
for all $\ell \in \{ 1,\ldots,L \} $, $g \in C^\infty_{\rm b}(U_\ell)$, 
$u \in L_q(U_\ell)$ and $v \in L_p(U_\ell)$.
Since $T$ is a pseudo-differential operator on the compact manifold $M$,
one can write 
\[
T = \sum_{\ell=1}^L M_{\chi_\ell} \, T \, M_{\widetilde \chi_\ell}
   + T_0 
,  \]
where $T_0$ has a $C^\infty$-kernel representation, i.e.,
there exists a $C^\infty$-function $K \colon M \times M \to \Ci$
such that 
\[
(T_0 u)(x)
= \int_M K(x,y) \, u(y) \, dy
\]
for all $u \in C_c^\infty(M)$ and $x \in M$.

The multi-commutator with $T_0$ is easy to estimate.
Let $g_1,\ldots,g_n \in W$.
Then 
\begin{eqnarray*}
|([M_{g_1},[\ldots,[M_{g_n}, T_0]\ldots]] u)(x)|
& = & \Big| \int_M K(x,y) \Big( \prod_{i=1}^n (g_i(x) - g_i(y)) \Big) \, u(y) \, dy \Big|  \\
& \leq & (3D)^n \int_M |K(x,y)| \, |u(y)| \, dy
\end{eqnarray*}
for all $u \in C^\infty(M)$ and $x \in M$, where we used 
Lemma~\ref{lpbdton300.5}\ref{lpbdton300.5-1}.
Hence
\[
\|[M_{g_1},[\ldots,[M_{g_n}, T_0]\ldots]] u\|_{L_q(M)}
\leq (3D)^n \, (\Vol(M))^{1 + \frac{1}{q} - \frac{1}{p}} \, \|K\|_\infty \, \|u\|_{L_p(M)}
\]
for all $u \in C^\infty(M)$.

Next we estimate the multi-commutators involving 
$M_{\chi_\ell} \, T \, M_{\widetilde \chi_\ell}$.
For all $\ell \in \{ 1,\ldots,L \} $ 
there exists a classical 
pseudo-differential operator $\widetilde T_\ell$ of order $k$ such that 
\[
\widetilde T_\ell w 
= \Bigg( \chi_\ell \, T \bigg( \Big( w \cdot (\widetilde \chi_\ell \circ \varphi_\ell^{-1}) \Big) 
           \circ \varphi_\ell \bigg) \Bigg) \circ \varphi_\ell^{-1}
\]
for all $w \in \cs(\Ri^m)$.
By the corresponding part of Theorem~\ref{tpbdton301}
there exists a $c_\ell > 0$ such that 
\[
\|[M_{h_1},[\ldots,[M_{h_n}, \widetilde T_\ell]\ldots]] u\|_q
\leq c_\ell \, \|\nabla h_1\|_\infty \ldots \|\nabla h_n\|_\infty \, \|u\|_p
\]
for all $h_1,\ldots,h_n \in \cs(\Ri^m)$ and $u \in C_c^\infty(\Ri^m)$.

Let $\ce \colon W^{1,\infty}(B(0,1)) \to W^{1,\infty}(\Ri^m)$ 
be an extension operator as in \cite{Ste1} Theorem~VI.5
with respect to the domain $B(0,1) \subset \Ri^m$.
Note that $\ce(h) \in C^\infty(\Ri^m)$ for all $h \in C^\infty(B(0,1))$.
Without loss of generality we may assume that 
$\supp \ce(h) \subset B(0,2)$ for all $h \in W^{1,\infty}(B(0,1))$.

Now let $\ell \in \{ 1,\ldots,L \} $.
Let $g \in W$.
Then 
\begin{eqnarray*}
\|\nabla \ce \Big( (g - g(\varphi_\ell^{-1}(0))) \circ \varphi_\ell^{-1} \Big)\|_{L_\infty(\Ri^m)}
& \leq & \|\ce\| \, \|g \circ \varphi_\ell^{-1} - (g \circ \varphi_\ell^{-1})(0)\|_{W^{1,\infty}(B(0,1))}  \\
& \leq & 2 \|\ce\| \, \|\nabla(g \circ \varphi_\ell^{-1})\|_{L_\infty(B(0,1))}  \\
& \leq & 2 c_0 \, \|\ce\| \, \|\nabla g\|_{L_\infty(M)}  \\
& \leq & 2 c_0 \, \|\ce\|
.
\end{eqnarray*}
Now let $g_1,\ldots,g_n \in W$.
Write $\check g_i = g_i - g_i(\varphi_\ell^{-1}(0))$
and $h_i = \ce(\check g_i \circ \varphi_\ell^{-1}) \in \cs(\Ri^m)$.
Then $\|\nabla h_i\|_{L_\infty(\Ri^m)} \leq 2 c_0 \, \|\ce\|$.
For all $A \subset \{ 1,\ldots,n \} $ define 
$\check g_A = \prod_{i \in A} \check g_i$ and $h_A = \prod_{i \in A} h_i$.
Let $u \in C^\infty(M)$.
Then 
\begin{eqnarray*}
[M_{g_1},[\ldots,[M_{g_n}, M_{\chi_\ell} \, T \, M_{\widetilde \chi_\ell}]\ldots]] u
& = & [M_{\check g_1},[\ldots,[M_{\check g_n}, M_{\chi_\ell} \, T \, M_{\widetilde \chi_\ell}]\ldots]] u \\
& = & \sum_{A \in \cp( \{ 1,\ldots,n \} )} (-1)^{n - |A|}
     \chi_\ell \, \check g_A \, T (\check g_{A^{\rm c}} \, \widetilde \chi_\ell \, u)
. 
\end{eqnarray*}
So 
\begin{eqnarray*}
\lefteqn{
\|[M_{g_1},[\ldots,[M_{g_n}, M_{\chi_\ell} \, T \, M_{\widetilde \chi_\ell}]\ldots]] u\|_{L_q(M)}
} \hspace*{30mm}  \\*
& \leq & c_0 \, \|\sum_{A \in \cp( \{ 1,\ldots,n \} )} (-1)^{n - |A|} \, 
    \Big( \chi_\ell \, \check g_A \, 
          T (\check g_{A^{\rm c}} \, \widetilde \chi_\ell \, u) 
                             \Big) \circ \varphi_\ell^{-1}\|_{L_q(\Ri^m)}  \\
& = & c_0 \, \|\sum_{A \in \cp( \{ 1,\ldots,n \} )} (-1)^{n - |A|} \, 
     (\check g_A \circ \varphi_\ell^{-1}) \, \widetilde T_\ell
         \Big( (\check g_{A^{\rm c}} \circ \varphi_\ell^{-1}) \cdot 
                 (u \circ \varphi_\ell^{-1}) \Big)\|_{L_q(\Ri^m)}  \\
& = & c_0 \, \|\sum_{A \in \cp( \{ 1,\ldots,n \} )} (-1)^{n - |A|} \, 
     h_A \, \widetilde T_\ell 
         \Big( h_{A^{\rm c}} \cdot (u \circ \varphi_\ell^{-1}) \Big)\|_{L_q(\Ri^m)}  \\
& = & c_0 \, 
   \|[M_{h_1},[\ldots,[M_{h_n},\widetilde T]\ldots]] (u \circ \varphi_\ell^{-1}) \|_{L_q(\Ri^m)}  \\
& \leq & c_0 \, c_\ell \, (2 c_0 \, \|\ce\|)^n \, \|(u \circ \varphi_\ell^{-1}) \|_{L_p(\Ri^m)}  \\
& \leq & c_0^2 \, c_\ell \, (2 c_0 \, \|\ce\|)^n \, \|u\|_{L_p(M)} 
.
\end{eqnarray*}
This proves the proposition.\hfill$\Box$

\vertspace

The proof of Theorem~\ref{tpbdton101} in the next section heavily depends on 
the bounds of the last proposition.

\section{Poisson bounds for $K^V_t$} \label{Spbdton4}

We assume for the rest of this paper that $\Omega \subset \Ri^d$ 
is bounded and connected, with a $C^\infty$-boundary $\Gamma$.
Recall that we do not assume that $\Gamma$ is connected. 
For the remaining part of this paper, fix an element in 
each connected component of $\Gamma$ as in Section~\ref{Spbdton3},
define $W$ as in (\ref{eSpbdton3;1}) and the distance 
$\rho_\Gamma$ as in (\ref{eSpbdton3;2}).
For all $g \in C^\infty(\Gamma)$ and $p \in [1,\infty]$ define the 
derivation $\delta_g$ on $\cl(L_p(\Gamma))$ by 
$\delta_g(E) = [M_g, E]$,where $M_g$ denotes the multiplication operator
with the function $g$.

In order not to repeat a proof for the kernel bound for $K_z$ with $z$ 
complex in Section~\ref{Spbdton6}, we prove a slightly more general 
proposition then that we need at the moment.
By Theorem~\ref{th2.5} we know that the assumptions of the 
next proposition are valid with $\alpha = 0$ and $N = 0$.
For all $\alpha \in [0,\frac{\pi}{2})$ define the sector
\begin{equation}
\Sigma_\alpha 
= \{ z \in \Ci : z = 0 \mbox{ or } |\arg z| \leq \alpha \} 
.  
\label{eSpbdton4;1}
\end{equation}
Note that $\Sigma_\alpha$ is closed.

\begin{prop} \label{ppbdton401}
For all $N \in [0,\infty)$ and $c > 0$ there exists a $c' > 0$ such that 
the following is valid.
Let $\alpha \in [0,\frac{\pi}{2})$ and suppose that 
\[
\|S_z\|_{p \to q}
\leq c \, (\cos \theta)^{-N} \, |z|^{-(d-1)(\frac{1}{p} - \frac{1}{q})}
\]
for all $p,q \in [1,\infty]$ and $z \in \Sigma_\alpha$, with $p \leq q$
and $0 < |z| \leq 1$, where $\theta = \arg z$.
Then 
\[
\|\delta_g^d(S_z)\|_{1 \to \infty}
\leq c' \, (\cos \theta)^{- N (d+1)} \, |z|
\]
for all $g \in W$ and $z \in \Sigma_\alpha$ with $0 < |z| \leq 1$,
where $\theta = \arg z$.
\end{prop}

For the proof we need the following decomposition for $\delta_g^d(S_z)$.
For all $k \in \Ni$ let 
\[
H_k = \{ (t_1,\ldots,t_{k+1}) \in (0,\infty)^{k+1} : t_1 + \ldots + t_{k+1} = 1 \}
\]
and let $d\lambda_k$ denote Lebesgue measure of the $k$-dimensional 
surface $H_k$.

\begin{lemma} \label{lpbdton402}
Let $T$ be a continuous semigroup on the sector $\Sigma_\alpha$
and generator $-A$ on a Banach space $\cx$, where $\alpha \in [0,\frac{\pi}{2})$.
Let $B \in \cl(\cx)$ and define the derivation $\delta$ on $\cl(\cx)$ by 
$\delta(E) = [B,E]$.
Then 
\begin{eqnarray*}
\delta^n(T_z)
& = & \sum_{k=1}^n (-z)^k
   \sum_{\scriptstyle j_1,\ldots,j_k \in \Ni \atop
         \scriptstyle j_1 + \ldots + j_k = n}
   \int_{H_k}
     T_{t_{k+1} \, z} \, \delta^{j_k}(A) \, T_{t_k \, z} 
       \circ \ldots \circ  \\*
& & \hspace*{50mm} {}
   \circ
     T_{t_2 \, z} \, \delta^{j_1}(A) \, T_{t_1 \, z} \, d\lambda_k(t_1,\ldots,t_{k+1})
\end{eqnarray*}
for all $z \in \Sigma_\alpha$ and $n \in \Ni$.
\end{lemma}
\proof\
If $n = 1$ then 
\[
\delta(T_z)
= [B,T_z]
= - z \int_0^1 T_{(1-s)z} \, [B,A] \, T_{sz} \, ds.
\]
Since $\delta$ is a derivation, the lemma easily follows by induction.\hfill$\Box$

\vertspace

\noindent
{\bf Proof of Proposition~\ref{ppbdton401}\hspace*{5pt}\ }
Recall that $\cn \in \OPS^1(M)$ (see \cite{Tay5} Appendix~C of Chapter~12).
By Proposition~\ref{ppbdton302} for all $p,q \in (1,\infty)$ with 
$p \leq q$ and $(d-1)(\frac{1}{p} - \frac{1}{q}) \in \{ 0,1,\ldots,d-1 \} $,
and in addition for the combination $p=1$ and $q=\infty$,
there exists a $c_{p,q} > 0$ such that 
\[
\|\delta_g^j(\cn)\|_{p \to q}
\leq c_{p,q}
\]
for all $g \in W$,
where $j = 1 + (d-1)(\frac{1}{p} - \frac{1}{q})$.

We will use the decomposition of Lemma~\ref{lpbdton402} and estimate
each term in the sum.
Let $k \in \{ 1,\ldots,d \} $, $(t_1,\ldots,t_{k+1}) \in H_k$,
$g \in W$ and $j_1,\ldots,j_k \in \Ni$ with 
$j_1 + \ldots + j_k = d$.

If $k = 1$ then $j_1 = d$ and 
\begin{eqnarray*}
|z|^k \, \|S_{t_2 z} \, \delta_g^{j_1}(\cn) \, S_{t_1 z}\|_{1 \to \infty}
& \leq & |z|^k \, \|S_{t_2 z}\|_{\infty \to \infty} \, 
      \|\delta_g^{j_1}(\cn)\|_{1 \to \infty} \, \|S_{t_1 z}\|_{1 \to 1}  \\
& \leq & c^2 \, c_{1,\infty} \, |z| \, (\cos \theta)^{-2N}
.  
\end{eqnarray*}
Suppose $k \in \{ 2,\ldots,d \} $.
There exists a $K \in \{ 1,\ldots,k+1 \} $ such that $t_K \geq \frac{1}{k+1}$.
Note that $\sum_{\ell = 1}^k (j_\ell - 1) = d-k < d-1$.
First suppose $K \not\in \{ 1,k+1 \} $.
Fix $1 = q_0 < p_1 \leq q_1 = p_2 \leq q_2 = p_3 \leq \ldots \leq 
q_{K-2} = p_{K-1} \leq q_{K-1} \leq p_K \leq q_K = p_{K+1} \leq q_{K+1} \leq \ldots \leq 
q_{k-1} = p_k \leq q_k < p_{k+1} = \infty$
such that 
\[
1 - \frac{1}{p_1} = \frac{1}{2(d-1)} = \frac{1}{q_k}
\quad , \quad
\frac{1}{p_\ell} - \frac{1}{q_\ell} = \frac{j_\ell - 1}{d-1}
\quad \mbox{and} \quad
\frac{1}{q_{K-1}} - \frac{1}{p_K} = \frac{k-2}{d-1}
\]
for all $\ell \in \{ 1,\ldots,k \} $.
Then 
\begin{eqnarray*}
\lefteqn{
|z|^k \, \|S_{t_{k+1} z} \, \delta_g^{j_k}(\cn)  \ldots \delta_g^{j_1}(\cn) \, S_{t_1 z}\|_{1 \to \infty}
} \hspace*{5mm}  \\*
& \leq & |z|^k \, \|S_{t_1 z}\|_{q_0 \to p_1}
       \prod_{\ell = 1}^k  \|S_{t_{\ell+1} z}\|_{q_\ell \to p_{\ell + 1}} 
                           \, \|\delta_g^{j_\ell}(\cn)\|_{p_\ell \to q_\ell} \\
& \leq & |z|^k \, c \, (\cos \theta)^{-N} \, (t_1 |z|)^{-(d-1)(\frac{1}{q_0} - \frac{1}{p_1})}
    \prod_{\ell = 1}^k c_{p_\ell , q_\ell} \, 
       c \, (\cos \theta)^{-N} \, (t_{\ell + 1} |z|)^{-(d-1)(\frac{1}{q_\ell} - \frac{1}{p_{\ell + 1}})}  \\
& = & c' \, (\cos \theta)^{-(k+1)N} \, |z|^k \, |z|^{-(k-1)} \, 
   t_1^{-1/2} \, t_K^{-(k-2)} \, t_{k+1}^{-1/2}  \\
& \leq & c' \, (k+1)^{k-2} \, (\cos \theta)^{-(k+1)N} \, |z| \, 
   t_1^{-1/2} \, t_{k+1}^{-1/2}
,
\end{eqnarray*}
where $c' = c^{k+1} \, \prod_{\ell = 1}^k c_{p_\ell , q_\ell}$.
If $K \in \{ 1,k+1 \} $ then a similar estimate is valid with possibly 
a different constant for $c'$.
Integration and taking the sum gives the proposition.\hfill$\Box$

\vertspace

We are now able to prove the Poisson bounds for real time.

\vertspace

\noindent
{\bf Proof of Theorem~\ref{tpbdton101}\hspace*{5pt}\ }
By Theorem~\ref{th2.5} and Proposition~\ref{ppbdton401}
there exists a $c > 0$ such that $\|\delta_g^d(S_t)\|_{1 \to \infty} \leq c \, t$
for all $g \in W$ and $t \in (0,1]$.
Hence 
\[
|(g(x) - g(y))^d \, K_t(x,y)| \leq c \, t
\]
for all $t \in (0,1]$, $x,y \in \Gamma$ and $g \in W$.
Optimising over $g \in W$ gives
$\rho_\Gamma(x,y)^d \, K_t(x,y) \leq c \, t$ and 
\[
\Big( \frac{\rho_\Gamma(x,y)}{t} \Big)^d K_t(x,y)
\leq c \, t^{-(d-1)}
\]
for all $x,y \in \Gamma$ and $t \in (0,1]$.
By Theorem~\ref{th2.5} there exists a $c_1 > 0$ such that 
$K_t(x,y) \leq \|S_t\|_{1 \to \infty} \leq c_1 \, t^{-(d-1)}$ for 
all $t \in (0,1]$ and $x,y \in \Gamma$.
Hence 
\[
\Big( 1 + \frac{\rho_\Gamma(x,y)}{t} \Big)^d K_t(x,y)
\leq 2^d \, (c_1 + c_2) \, t^{-(d-1)}
.  \]
Since $\rho_\Gamma$ is equivalent to the distance $(x,y) \mapsto |x-y|$
on $\Gamma$ by the Lemma~\ref{lpbdton300.5}\ref{lpbdton300.5-4}, 
one establishes that there is a
$c_2 > 0$ such that 
\[
K^V_t(x,y) \le K_t(x,y)
\leq \frac{c_2 \, t^{-(d-1)} }
          {\displaystyle \Big( 1 + \frac{|x-y|}{t} \Big)^d }
\]
for all $t \in (0,1]$ and $x,y \in \Gamma$, where 
we used the domination of Theorem~\ref{th2.3} in the first inequality.

Finally we deduce large time bounds.
Using Theorem~\ref{th2.5} there is a $c_3 > 0$ such that 
\[
\|S^V_t\|_{1 \to \infty} \leq c_3 \, (t \wedge 1)^{-(d-1)} \, e^{-\lambda_1 t}
\]
for all $t \in [1,\infty)$.
Since $\Gamma$ is bounded, there is a $c_4 > 0$ such that 
\[
K^V_t(x,y)
\leq \frac{c_4 \, (t \wedge 1)^{-(d-1)} \, e^{-\lambda_1 t}}
          {\displaystyle \Big( 1 + \frac{|x-y|}{t} \Big)^d }
\]
for all $x,y \in \Gamma$ and $t \in [1,\infty)$.
This completes the proof of Theorem~\ref{tpbdton101}.\hfill$\Box$

\section{Poisson bounds for $K_z$} \label{Spbdton6}

In this section we will give a proof for Theorem~\ref{tpbdton102}, that is 
Poisson kernel bounds for complex time.
The proof follows from Proposition~\ref{ppbdton401}, once one has
semigroup bounds for $\|S_z\|_{p \to q}$ for all $1 \leq p \leq q \leq \infty$.
These bounds are easy if $p \leq 2 \leq q$, see Lemma~\ref{lpbdton602}.
But if $2 \not\in [p,q]$ then it is much harder.
The  method to derive them is to prove bounds for $\delta^d_g(S_z)$
from $L_1$ to $C^\nu = W^{1,p}$. 
Unfortunately, this method does not allow to give directly the 
bounds from $L_1$ to $L_\infty$.
It is convenient to consider the semigroup generated by a power of 
$\cn$ and then use fractional powers to go back to $\cn$.

Define $P = \cn + I$.
If confusion is possible, then we write $P_p$ for the operator on 
$L_p(\Gamma)$, where $p \in [1,\infty]$.
We start with a regularity result for the Dirichlet-to-Neumann operator.

\begin{prop} \label{ppbdton601}
Let $p \in (1,\infty)$ and $n \in \Ni_0$.
Then $W^{n,p}(\Gamma) = D(P_p^n)$.
In particular, there exists a $c > 0$ such that 
\[
c^{-1} \, \|u\|_{W^{n,p}(\Gamma)} 
\leq \|P_p^n u\|_p
\leq c \, \|u\|_{W^{n,p}(\Gamma)} 
\]
for all $u \in W^{n,p}(\Gamma)$.
\end{prop}
\proof\
The case $n = 0$ is trivial.
Let $n \in \Ni_0$ and suppose that $W^{n,p}(\Gamma) = D(P_p^n)$.
It follows from (C.4) or Proposition~C.1 in Appendix~C of Chapter~12
in \cite{Tay5} that there exists a pseudo-differential operator $V_0$
of order $0$ such that $P = \sqrt{-\Delta} + V_0$.
Then $P^{n+1} = (-\Delta)^{(n+1)/2} + W$, where $W \in \OPS^n(\Gamma)$.
By Lemma~\ref{lpbdtonA03} there exists a $c > 0$ such that 
$\|W u\|_p \leq c \, \|u\|_{W^{n,p}(\Gamma)}$ for all $u \in C^\infty(\Gamma)$.
By Proposition~\ref{ppbdtonA01} one has $W^{n+1,p}(\Gamma) = D((-\Delta_p)^{(n+1)/2})$
with equivalent norms.
Hence there exists a $c' > 0$ such that 
$\|(-\Delta_p)^{(n+1)/2} u\|_p \leq c' \, \|u\|_{W^{n+1,p}(\Gamma)}$
for all $u \in C^\infty(\Gamma)$.
Then $\|P^{n+1} u\|_p \leq (c+c') \, \|u\|_{W^{n+1,p}(\Gamma)}$
for all $u \in C^\infty(\Gamma)$.
Since $C^\infty(\Gamma)$ is dense in $W^{n+1,p}(\Gamma)$ (see Lemma~\ref{lpbdtonA03})
and $P$ is closed, 
it follows that $W^{n+1,p}(\Gamma) \subset D(P_p^{n+1})$.
The converse follows similarly, once one knows that $C^\infty(\Gamma)$ is  
a core for $P_p^{n+1}$.
The latter can be proved as follows.
Let $m \in \Ni$.
Then $P_2^m$ is an elliptic pseudo-differential operator of order~$m$.
Hence $D(P_2^m) = W^{m,2}(\Gamma)$ by \cite{Kum} Theorem~3.6.7.
So if $S^{(p)}$ denotes the semigroup generated by $-P_p^{n+1}$,
then 
\[
S^{(p)}(C^\infty(\Gamma))
= S^{(2)}(C^\infty(\Gamma))
\subset \bigcap_{m=1}^\infty D((P_2^{n+1})^m)
= \bigcap_{m=1}^\infty W^{(n+1)m,2}(\Gamma)
= C^\infty(\Gamma)
,  \]
where we used the Sobolev embedding of Proposition~\ref{ppbdtonA02}
in the last step.
Hence $C^\infty(\Gamma)$ is a core for $P_p^{n+1}$ and the proof of the
proposition is complete.\hfill$\Box$

\vertspace

Let $S$ be the semigroup generated by $-P$.
For all $m \in \Ni$
let $T^{(m)}$ be the semigroup on $L_2(\Gamma)$ generated by $-P^m = - (\cn + I)^m$.
Clearly $T^{(m)}$ is holomorphic with angle $\pi/2$.

\begin{lemma} \label{lpbdton602}
Let $m\in \Ni$, $n \in \Ni_0$ and $p \in (2,\infty]$.
Then there exists a 
$c > 0$ such that $T^{(m)}_z (L_2(\Gamma)) \subset C^\infty(\Gamma)$ and 
\[
\|T^{(m)}_z\|_{L_2 \to W^{n,p}}
\leq c \, |\RRe z|^{-\frac{d-1}{m}(\frac{1}{2} - \frac{1}{p})} \, |\RRe z|^{-\frac{n}{m}} 
\]
for all $z \in \Ci$ with $\RRe z > 0$.
\end{lemma}
\proof\
Clearly 
$T^{(m)}_z (L_2(\Gamma))
\subset \bigcap_{\ell=1}^\infty D(P^{m \ell})
= \bigcap_{\ell=1}^\infty W^{m \ell,2}(\Gamma)
= C^\infty(\Gamma)$ by Proposition~\ref{ppbdton601} and the Sobolev embedding 
of Proposition~\ref{ppbdtonA02}.
In addition, $D(P^{d-1+n}) = W^{d-1+n,2}(\Gamma) \subset W^{n,p}(\Gamma)$.
By Propositions~\ref{ppbdtonA02} and \ref{ppbdton601} there exists a $c > 0$
such that 
\[
\|u\|_{W^{n,p}} 
\leq c \, \|P^{d-1+n} u\|_2^\alpha \, 
          \|u\|_2^{1-\alpha}
\]
for all $u \in C^\infty(\Gamma)$, where 
$\alpha = \dfrac{n + (d-1)(\frac{1}{2} - \frac{1}{p})}{n+d-1}$.
Then the lemma follows by the spectral theorem.\hfill$\Box$

\vertspace

We will use again Lemma~\ref{lpbdton402} to decompose $\delta_g^d(T^{(m)}_z)$.
This time it involves higher order derivatives on $g$.
For all $k \in \Ni$ define $W_k$ as in (\ref{eSpbdton3;3}).
In order to estimate $\delta_g^j(P^m) \, T^{(m)}_z$ we need a few lemmas.
The third one is the most delicate.

\begin{lemma} \label{lpbdton602.3}
Let  $\alpha$ be a multi-index over 
$ \{ 1,\ldots,d-1 \} $ and let $j \in \Ni$ with $|\alpha| \leq j$.
Then there exist constants $c_{\alpha_1,\ldots,\alpha_{k+1}} \in \Ri$, 
where $k \in \{ 0,\ldots,|\alpha| \} $ and 
$\alpha_1,\ldots,\alpha_{k+1}$ are multi-indices,
such that 
\[
\partial^\alpha \, \delta_h^j(T)
= \sum_{k=0}^{|\alpha|}
  \sum_{{\scriptstyle \alpha_1,\ldots,\alpha_{k+1} \atop
         \scriptstyle |\alpha_1|,\ldots,|\alpha_k| \geq 1} \atop
         \scriptstyle |\alpha_1| + \ldots + |\alpha_{k+1}| = |\alpha| }
    c_{\alpha_1,\ldots,\alpha_{k+1}} \, 
     M_{\partial^{\alpha_1} h} \ldots M_{\partial^{\alpha_k} h}
        \delta_h^{j-k}(\partial^{\alpha_{k+1}} T)
\]
for every $h \in \cs(\Ri^{d-1})$ and pseudo-differential operator $T$.
\end{lemma}
\proof\
It follows by induction to $j$ that 
$\partial_i \, \delta_h^j(T) 
= j \, M_{\partial_i h} \, \delta_h^{j-1}(T) + \delta_h^j(\partial_i T)$
for all $i \in \{ 1,\ldots,d-1 \} $ and $j \in \Ni$.
Then the lemma follows by induction to $|\alpha|$.\hfill$\Box$

\vertspace

In the next lemma we move the derivatives to the right.

\begin{lemma} \label{lpbdton602.5}
Let $j \in \Ni$ and let $\beta$ be a multi-index over 
$ \{ 1,\ldots,d-1 \} $.
Then there exist constants $\tilde c_{\beta,\ldots,\beta_{j+2}} \in \Ri$,
where $\beta,\ldots,\beta_{j+2}$ are multi-indices,
such that 
\[
\partial^\beta \, \delta_h^j(T)
= \sum_{\scriptstyle \beta_1,\ldots,\beta_{j+2} \atop
        \scriptstyle |\beta_1| + \ldots + |\beta_{j+2}| = |\beta|}
    \tilde c_{\beta_1,\ldots,\beta_{j+2}} \, 
     (\delta_{\partial^{\beta_1} h} \ldots 
           \delta_{\partial^{\beta_j} h}(\partial_{\beta_{j+1}} T))
    \circ \partial^{\beta_{j+2}}
\]
for every $h \in \cs(\Ri^{d-1})$ and pseudo-differential operator $T$, 
where 
\[
\partial_{\beta_{j+1}} T = [\partial_{i_1},[\ldots,[\partial_{i_k},T] \ldots]]
\]
if $\beta_{j+1} = (i_1,\ldots,i_k)$.
\end{lemma}
\proof\
Since $[\partial_i ,\delta_h(T)] = \delta_{\partial_i h}(T) + \delta_h([\partial_i,T])$,
the lemma easily follows by induction to $j$ and $|\beta|$.\hfill$\Box$

\vertspace

The next lemma is the key estimate in our proof to estimate
$\|\delta^d(T^{(m)}_z)\|_{L_1 \to C^\nu}$.

\begin{lemma} \label{lpbdton603}
For all $m_1,m_2 \in \Ni_0$ and $j \in \Ni$ with $m_1 + m_2 + 1 \geq j$
there exists a $c > 0$ such that 
\begin{equation}
\|P^{m_1} \, \delta_g^j(P) \, P^{m_2} u\|_2
\leq c \, \|P^{m_1 + m_2 + 1 - j} u\|_2
\label{elpbdton603;11}
\end{equation}
for all $u \in C^\infty(\Gamma)$ and $g \in W_{m_1 + m_2 + 1}$.
\end{lemma}
\proof\
We may assume that $m_2 = 0$, or $m_1 + m_2 + 1 = j$.

We use the notation as in the proof of Proposition~\ref{ppbdton302} with $p = q = 2$
and with $T = P$.
Now $m = d-1$.
Thus let $L \in \Ni$, $c_0 > 0$, $T_0$, $K \in C^\infty(\Gamma \times \Gamma)$
and for all $\ell \in \{ 1,\ldots,L \} $ let $U_\ell$, $\varphi_\ell$,
$\chi_\ell$, $\widetilde \chi_\ell$ and
$\widetilde T_\ell$ be as in the proof of 
Proposition~\ref{ppbdton302}.
We may assume that 
\[
\sum_{i=1}^{m_1 + m_2 + 1} \|\nabla^i(g \circ \varphi_\ell^{-1})\|_{L_\infty(B(0,1))}
\leq c_0 \sum_{i=1}^{m_1 + m_2 + 1} \|\nabla^i g\|_{L_\infty(\Gamma)}
\]
for all $\ell \in \{ 1,\ldots,L \} $ and $g \in C^\infty_{\rm b}(U_\ell)$.
Moreover, let $\hat \chi_\ell \in C_c^\infty(U_\ell)$ be such that 
$\hat \chi_\ell(x) = 1$ for all $x \in \supp \hat \chi_\ell$.
Then 
\[
P = \sum_{\ell=1}^L M_{\chi_\ell} \, P \, M_{\widetilde \chi_\ell}
   + T_0 
,  \]
where $T_0$ has $K$ as kernel.

We first estimate the contribution of the operator $T_0$ in (\ref{elpbdton603;11}).
Note that 
\[
\|P^{m_1} \, \delta_g^j(T_0) \, P^{m_2} u\|_2
\leq \sup_{v \in C^\infty(\Gamma), \; \|v\|_2 \leq 1}
   \sum_{\ell_1 = 1}^L
   \sum_{\ell_2 = 1}^L
     |(\delta^j_g(T_0) \, M_{\chi_{\ell_2}} \, P^{m_2} u, M_{\chi_{\ell_1}} \, P^{m_1} v)|
.  \]
Let $\ell_1,\ell_2 \in \{ 1,\ldots,L \} $.
By Lemma~\ref{lpbdtonA03.5} for every multi-index $\gamma$
over $ \{ 1,\ldots,d-1 \} $ with $|\gamma| \leq m_1$ there exists 
a bounded operator $T^{(1)}_\gamma$ on $L_2(\Gamma)$ such that 
\[
M_{\chi_{\ell_1}} \, P^{m_1} 
= \sum_{|\gamma| \leq m_1} 
   M_{\chi_{\ell_1}} \, \Big( \frac{\partial}{\partial \varphi_{\ell_1}} \Big)^\gamma \, T^{(1)}_\gamma
.  \]
Similarly write 
\[
M_{\chi_{\ell_2}} \, P^{m_2} 
= \sum_{|\gamma| \leq m_2} 
   M_{\chi_{\ell_2}} \, \Big( \frac{\partial}{\partial \varphi_{\ell_2}} \Big)^\gamma \, T^{(2)}_\gamma
\]
with $T^{(2)}_\gamma \in \cl(L_2(\Gamma))$.
By (\ref{eSpbdtonapp;1}) there exists a $c_1 \geq 1$ such that 
\begin{equation}
\Big\| \one_{U_{\ell_1}} \, \Big( \frac{\partial}{\partial \varphi_{\ell_1}} \Big)^\gamma g \Big\|_\infty
\leq c_1 \, \|\nabla^{|\gamma|} g\|_\infty
\quad \mbox{and} \quad
\Big\| \one_{U_{\ell_2}} \, \Big( \frac{\partial}{\partial \varphi_{\ell_2}} \Big)^\gamma g \Big\|_\infty
\leq c_1 \, \|\nabla^{|\gamma|} g\|_\infty
\label{elpbdton603;12}
\end{equation}
for all $g \in W^{m_1 + m_2 + 1,\infty}(\Gamma)$ and $|\gamma| \leq m_1 + m_2 + 1$.
Let $u,v \in C^\infty(\Gamma)$ and $g \in W_{m_1 + m_2 + 1}$.
Then 
\begin{eqnarray*}
\lefteqn{
|(\delta^j_g(T_0) \, M_{\chi_{\ell_2}} \, P^{m_2} u, M_{\chi_{\ell_1}} \, P^{m_1} v)|
} \hspace{30mm} \\*
& \leq & \sum_{|\gamma_1| \leq m_1}
         \sum_{|\gamma_2| \leq m_2}
   |( \Big( \frac{\partial}{\partial \varphi_{\ell_1}} \Big)^{\gamma_1} \, 
      M_{\chi_{\ell_1}} \, \delta^j_g(T_0) \, M_{\chi_{\ell_2}} \, 
       \Big( \frac{\partial}{\partial \varphi_{\ell_2}} \Big)^{\gamma_2} \, 
      T^{(2)}_{\gamma_2} u, T^{(1)}_{\gamma_1} v)|
.  
\end{eqnarray*}
Note that 
\[
(M_{\chi_{\ell_1}} \, \delta^j_g(T_0) \, M_{\chi_{\ell_2}} w)(x)
= \int_\Gamma \chi_{\ell_1}(x) \, (g(x) - g(y))^j \, K(x,y) \, \chi_{\ell_2}(y) \, w(y) \, dy
\]
for all $x \in \Gamma$ and $w \in C^\infty(\Gamma)$.
Moreover, $|g(x) - g(y)| \leq 3D$ for all $x,y \in \Gamma$ by 
Lemma~\ref{lpbdton300.5}\ref{lpbdton300.5-1}.
Using (\ref{elpbdton603;12}) and the product rule one estimates
\begin{eqnarray*}
\lefteqn{
\| \Big( \frac{\partial}{\partial \varphi_{\ell_1}} \Big)^{\gamma_1} \, 
      M_{\chi_{\ell_1}} \, \delta^j_g(T_0) \, M_{\chi_{\ell_2}} \, 
       \Big( \frac{\partial}{\partial \varphi_{\ell_2}} \Big)^{\gamma_2} w\|_2
} \hspace{1mm} \\*
& \leq & (c_1 (j+2))^{|\gamma_1| + |\gamma_2|} (1 + 3D)^j
   \|\chi_{\ell_1}\|_{W^{|\gamma_1|,\infty}} \, \|\chi_{\ell_2}\|_{W^{|\gamma_2|,\infty}} \,
     \Vol(\Gamma) \sum_{i=0}^{|\gamma_1|} \sum_{i'=0}^{|\gamma_1|}
        \|\nabla_{(1)}^i \nabla_{(2)}^{i'} K\|_\infty \, \|w\|_2
\end{eqnarray*}
for all $w \in C^\infty(\Gamma)$.
Now it is clear that there exists a $c_2 > 0$ such that 
\[
|(\delta^j_g(T_0) \, M_{\chi_{\ell_1}} \, P^{m_2} u, M_{\chi_{\ell_2}} \, P^{m_1} v)|
\leq c_2 \, \|u\|_2 \, \|v\|_2
\]
for all $u,v \in C^\infty(\Gamma)$ and $g \in W_{m_1 + m_2 + 1}$.
Then 
\[
\|P^{m_1} \, \delta_g^j(T_0) \, P^{m_2} u\|_2
\leq c_2 \, L^2 \, \|u\|_2
\leq c_2 \, L^2 \, \|P^{m_1 + m_2 + 1 - j} u\|_2
\]
for all $u \in C^\infty(\Gamma)$ and $g \in W_{m_1 + m_2 + 1}$.

The estimates for the other terms in the decomposition of $P$ 
involve much more work, as in Proposition~\ref{ppbdton302}.
This time
let $\ce \colon W^{m_1 + m_2 + 1,\infty}(B(0,1)) \to W^{m_1 + m_2 + 1,\infty}(\Ri^{d-1})$ 
be an extension operator as in \cite{Ste1} Theorem~VI.5
with respect to the domain $B(0,1) \subset \Ri^{d-1}$.
Again note that $\ce(h) \in C^\infty(\Ri^{d-1})$ for all $h \in C^\infty(B(0,1))$.
Without loss of generality we may assume that 
$\supp \ce(h) \subset B(0,2)$ for all $h \in W^{m_1 + m_2 + 1,\infty}(B(0,1))$.
Let $\ell \in \{ 1,\ldots,L \} $.
Let $g \in W_{m_1 + m_2 + 1}$.
Then 
\begin{eqnarray*}
\|\nabla^i \ce \Big( (g - g(\varphi_\ell(0))) \circ \varphi_\ell^{-1} \Big)\|_{L_\infty(\Ri^{d-1})}
& \leq & \|\ce\| \, \|g \circ \varphi_\ell^{-1} - (g \circ \varphi_\ell^{-1})(0)\|_{W^{m_1 + m_2 + 1,\infty}(B(0,1))}  \\
& \leq & 2 \|\ce\| \, \sum_{i'=1}^{m_1 + m_2 + 1} \|\nabla^{i'} (g \circ \varphi_\ell^{-1})\|_{L_\infty(B(0,1))}  \\
& \leq & 2 c_0 \, \|\ce\| \, \sum_{i'=1}^{m_1 + m_2 + 1} \|\nabla^{i'} g\|_{L_\infty(\Gamma)}  \\
& \leq & C
\end{eqnarray*}
for all $i \in \{ 1,\ldots,m_1 + m_2 + 1 \} $, where 
$C = 2 c_0 \, \|\ce\| \, (m_1 + m_2 + 1)$.

As a consequence of Proposition~\ref{ppbdton601} there exists a $c > 0$
such that 
$\|P^{m_1}  u\|_2 
    \leq c \, \sum_{i=0}^{m_1} \|\nabla^i  u\|_2$
for all $u \in C^\infty(\Gamma)$.
So it suffices to show that there exists a $c > 0$ such that 
\begin{equation}
\|\nabla^i \, \delta_g^j(M_{\chi_\ell} \, P \, M_{\widetilde \chi_\ell}) \, P^{m_2} u\|_2
\leq c \, \|P^{m_1 + m_2 + 1 - j} u\|_2
\label{elpbdton603;10}
\end{equation}
for all $u \in C^\infty(\Gamma)$, $g \in W_{m_1 + m_2 + 1}$, $\ell \in \{ 1,\ldots,L \} $
and $i \in \{ 0,\ldots,m_1 \} $.
Next fix $\ell \in \{ 1,\ldots,L \} $.

First suppose that $m_2 = 0$.
Let $\alpha,\beta$ be a multi-indices over $ \{ 1,\ldots,d-1 \} $ with 
$|\alpha| \leq j-1$ and $|\beta| \leq (m_1 + 1 - j) \vee 0$.
Let $g \in W_{m_1 + m_2 + 1}$.
Choose $h = \ce(\check g \circ \varphi_\ell^{-1})$ where $\check g = g - g(\varphi_\ell^{-1}(0))$.
Using Lemma~\ref{lpbdton602.3} one has
\begin{eqnarray*}
\lefteqn{
\|\partial_{\varphi_\ell}^\beta \, \partial_{\varphi_\ell}^\alpha 
    \delta_g^j(M_{\chi_\ell} \, P \, M_{\widetilde \chi_\ell}) u\|_{L_2(\Gamma)}
} \hspace{2mm}  \\
& \leq & c_0 \, \|(\partial^\beta \, \partial^\alpha \, \delta_h^j(\widetilde T_\ell)) (u \circ \varphi_\ell^{-1})\|_{L_2(\Ri^{d-1})}  \\
& \leq & c_0 \sum_{k=0}^{|\alpha|} \hspace{-2mm}
  \sum_{{\scriptstyle \alpha_1,\ldots,\alpha_{k+1} \atop
         \scriptstyle |\alpha_1|,\ldots,|\alpha_k| \geq 1} \atop
         \scriptstyle |\alpha_1| + \ldots + |\alpha_{k+1}| = |\alpha| }
    |c_{\alpha_1,\ldots,\alpha_{k+1}}| \, 
     \|(\partial^\beta M_{\partial^{\alpha_1} h} \ldots M_{\partial^{\alpha_k} h}
        \delta_h^{j-k}(\partial^{\alpha_{k+1}} \widetilde T_\ell))
           (u \circ \varphi_\ell^{-1})\|_{L_2(\Ri^{d-1})}  \\
& \leq & c_0 \sum_{|\gamma| \leq |\beta|}
  \sum_{k=0}^{|\alpha|} \hspace{-2mm}
  \sum_{{\scriptstyle \alpha_1,\ldots,\alpha_{k+1} \atop
         \scriptstyle |\alpha_1|,\ldots,|\alpha_k| \geq 1} \atop
         \scriptstyle |\alpha_1| + \ldots + |\alpha_{k+1}| = |\alpha| }
    |c_{\alpha_1,\ldots,\alpha_{k+1}}| \, 
      (k+1)^{|\beta|} \, C^k 
     \|(\partial^\gamma 
        \delta_h^{j-k}(\partial^{\alpha_{k+1}} \widetilde T_\ell))
           (u \circ \varphi_\ell^{-1})\|_{L_2(\Ri^{d-1})}  
.
\end{eqnarray*}
But then Lemma~\ref{lpbdton602.5} gives
\begin{eqnarray*}
\lefteqn{
\|(\partial^\gamma 
        \delta_h^{j-k}(\partial^{\alpha_{k+1}} \widetilde T_\ell))
           (u \circ \varphi_\ell^{-1})\|_{L_2(\Ri^{d-1})}  
} \hspace{5mm}  \\
& \leq & 
     \sum_{\scriptstyle \beta_1,\ldots,\beta_{j-k+2} \atop
           \scriptstyle |\beta_1| + \ldots + |\beta_{j-k+2}| = |\gamma|}
        |\tilde c_{\beta_1,\ldots,\beta_{j-k+2}}| \cdot  \\*
& & \hspace{30mm} {} \cdot
     \|((\delta_{\partial^{\beta_1} h} \ldots 
       \delta_{\partial^{\beta_{j-k}} h}
          (\partial_{\beta_{j-k+1}}(\partial^{\alpha_{k+1}} \, \widetilde T_\ell) ))
               \circ \partial^{\beta_{j-k+2}}) (u \circ \varphi_\ell^{-1})\|_{L_2(\Ri^{d-1})}   \\
& \leq & c_3 
     \sum_{\scriptstyle \beta_1,\ldots,\beta_{j-k+2} \atop
        \scriptstyle |\beta_1| + \ldots + |\beta_{j-k+2}| = |\gamma|}
        |\tilde c_{\beta_1,\ldots,\beta_{j-k+2}}| \:
   \|\nabla \partial^{\beta_1} h\|_\infty \ldots 
   \|\nabla \partial^{\beta_{j-k}} h\|_\infty  \cdot  \\*
& & \hspace{30mm} {} \cdot
       \|\partial^{\beta_{j-k+2}} ((\hat \chi_\ell \, u) \circ \varphi_\ell^{-1})\|_{L_2(\Ri^{d-1})}  \\
& \leq & C^{j-k} \, c_3 
     \sum_{\scriptstyle \beta_1,\ldots,\beta_{j-k+2} \atop
        \scriptstyle |\beta_1| + \ldots + |\beta_{j-k+2}| = |\gamma|}
        |\tilde c_{\beta_1,\ldots,\beta_{j+2}}| \, 
\|\partial^{\beta_{j-k+2}} ((\hat \chi_\ell \, u) \circ \varphi_\ell^{-1})\|_{L_2(\Ri^{d-1})}
\end{eqnarray*}
for a suitable $c_3 > 0$,
where we used the Coifman--Meyer estimate of Theorem~\ref{tpbdton301}\ref{tpbdton301-1}
in the penultimate step.
This is possible since $|\alpha_{k+1}| \leq |\alpha| - k \leq j-k-1$ and hence
$\partial^{\alpha_{k+1}} \widetilde T_\ell \in OPS^{j-k}$
and then also $\partial_{\beta_{j-k+1}}(\partial^{\alpha_{k+1}} \, \widetilde T_\ell) \in \OPS^{j-k}$
by \cite{Ste3} Theorem~VI.7.3.
Then $\|\partial^{\beta_{j-k+2}} ((\hat \chi_\ell \, u) \circ \varphi_\ell^{-1})\|_{L_2(\Ri^{d-1})} 
\leq c_4 \, \|P^{|\beta_{j-k+2}} u\|_{L_2(\Gamma)} \leq c_4 \, \|P^{m_1 + 1 - j} u\|_{L_2(\Gamma)}$
for a suitable $c_4 > 0$.
This completes the proof of (\ref{elpbdton603;10}) if $m_2 = 0$.

Finally suppose that $m_1 + m_2 + 1 = j$.
Note that $P^{m_2} \in \OPS^{m_2}$.
Using Lemma~\ref{lpbdtonA03.5}
it follows that for every multi-index $\gamma$ over $ \{ 1,\ldots,d-1 \} $
with $|\gamma| \leq m_2$ there exists a bounded operator $T_\gamma$ 
on $L_2(\Gamma)$ such that 
\[
M_{\widetilde \chi_\ell} \, P^{m_2}
= \sum_{|\gamma| \leq m_2}
   M_{\widetilde \chi_\ell} \, \Big( \frac{\partial}{\partial \varphi_\ell} \Big)^\gamma \, T_\gamma
.  \]
Let $\alpha$ be a multi-index with $|\alpha| \leq m_1$.
We shall show that (\ref{elpbdton603;10}) is valid.
Using Lemma~\ref{lpbdton602.3} twice one deduces that first
\begin{eqnarray*}
\lefteqn{
\|\partial_{\varphi_\ell}^\alpha \, 
      \delta_g^j(M_{\chi_\ell} \, P \, M_{\widetilde \chi_\ell}) \, P^{m_2} u)\|_{L_2(\Gamma)}
} \hspace{5mm}  \\*
& \leq & c_0 \sum_{|\gamma| \leq m_2}
   \|\partial^\alpha \, \delta_h^j(\widetilde T_\ell) \, \partial^\gamma ((T_\gamma u) \circ \varphi_\ell^{-1})\|_{L_2(\Ri^{d-1})}  \\
& \leq & c_0 \sum_{|\gamma| \leq m_2}
             \sum_{k=0}^{|\alpha|}
  \sum_{{\scriptstyle \alpha_1,\ldots,\alpha_{k+1} \atop
         \scriptstyle |\alpha_1|,\ldots,|\alpha_k| \geq 1} \atop
         \scriptstyle |\alpha_1| + \ldots + |\alpha_{k+1}| = |\alpha| }
    |c_{\alpha_1,\ldots,\alpha_{k+1}}| \cdot  \\*
& & \hspace{50mm} {} \cdot
     \|M_{\partial^{\alpha_1} h} \ldots M_{\partial^{\alpha_k} h}
        \delta_h^{j-k}(\partial^{\alpha_{k+1}} \widetilde T_\ell) \, 
      \partial^\gamma ((T_\gamma u) \circ \varphi_\ell^{-1})\|_{L_2(\Ri^{d-1})}   \\
& \leq & c_0 \sum_{|\gamma| \leq m_2}
             \sum_{k=0}^{|\alpha|}
  \sum_{{\scriptstyle \alpha_1,\ldots,\alpha_{k+1} \atop
         \scriptstyle |\alpha_1|,\ldots,|\alpha_k| \geq 1} \atop
         \scriptstyle |\alpha_1| + \ldots + |\alpha_{k+1}| = |\alpha| }
    |c_{\alpha_1,\ldots,\alpha_{k+1}}| \, C^k \, 
        \|\delta_h^{j-k}(\partial^{\alpha_{k+1}} \widetilde T_\ell) \, 
      \partial^\gamma ((T_\gamma u) \circ \varphi_\ell^{-1})\|_{L_2(\Ri^{d-1})}  
\end{eqnarray*}
and next 
\begin{eqnarray*}
\lefteqn{
\|\delta_h^{j-k}(\partial^{\alpha_{k+1}} \widetilde T_\ell) \, 
      \partial^\gamma ((T_\gamma u) \circ \varphi_\ell^{-1})\|_{L_2(\Ri^{d-1})}  
} \hspace{5mm}  \\*
& \leq & \sum_{k'=0}^{|\gamma|}
  \sum_{{\scriptstyle \gamma_1,\ldots,\gamma_{k'+1} \atop
         \scriptstyle |\gamma_1|,\ldots,|\gamma_{k'}| \geq 1} \atop
         \scriptstyle |\gamma_1| + \ldots + |\gamma_{k'+1}| = |\gamma| }
    |c_{\gamma_1,\ldots,\gamma_{k'+1}}|  \cdot  \\*
& & \hspace{25mm} {} \cdot
  \|\delta_h^{j-k-k'}(\partial^{\alpha_{k+1}} \widetilde T_\ell \partial^{\gamma_{k'+1}} ) \, 
   M_{\partial^{\gamma_{k'}} h} \ldots M_{\partial^{\gamma_1} h} \, 
    ((T_\gamma u) \circ \varphi_\ell^{-1})\|_{L_2(\Ri^{d-1})}    \\
& \leq & c_5 \sum_{k'=0}^{|\gamma|}
  \sum_{{\scriptstyle \gamma_1,\ldots,\gamma_{k'+1} \atop
         \scriptstyle |\gamma_1|,\ldots,|\gamma_{k'}| \geq 1} \atop
         \scriptstyle |\gamma_1| + \ldots + |\gamma_{k'+1}| = |\gamma| }
    |c_{\gamma_1,\ldots,\gamma_{k'+1}}| \, C^{k'} \, \|\nabla h\|_\infty^{j-k-k'} \, 
  \|((T_\gamma u) \circ \varphi_\ell^{-1})\|_{L_2(\Ri^{d-1})}   \\
& \leq & c_0 \, c_5 \sum_{k'=0}^{|\gamma|}
  \sum_{{\scriptstyle \gamma_1,\ldots,\gamma_{k'+1} \atop
         \scriptstyle |\gamma_1|,\ldots,|\gamma_{k'}| \geq 1} \atop
         \scriptstyle |\gamma_1| + \ldots + |\gamma_{k'+1}| = |\gamma| }
    |c_{\gamma_1,\ldots,\gamma_{k'+1}}| \, C^{j-k} \, 
    \|T_\gamma\|_{2 \to 2} \, 
     \|u\|_{L_2(\Gamma)} 
\end{eqnarray*}
for a suitable $c_5 > 0$,
where we used again the Coifman--Meyer estimate of Theorem~\ref{tpbdton301}\ref{tpbdton301-1}
in the penultimate step.
This is possible since 
$|\alpha_{k+1}| + 1 + |\gamma_{k'+1}| \leq m_1 - k + 1 + m_2 - k' = j - k - k'$
and hence 
$\partial^{\alpha_{k+1}} \, \widetilde T_\ell \, \partial^{\gamma{k'+1}} \in OPS^{j-k-k'}$.
The proof of the Lemma~\ref{lpbdton603} is complete.\hspace*{10mm}\hfill$\Box$

\begin{lemma} \label{lpbdton604}
For all $n \in \Ni$, $m_1,\ldots,m_{n+1} \in \Ni_0$ and $j_1,\ldots,j_n \in \Ni$
with $m_1 + \ldots + m_{n+1} + n \geq j_1 + \ldots + j_n$ there exists a $c > 0$
such that 
\[
\|P^{m_1} \, \delta_g^{j_1}(P) \, P^{m_2} \ldots P^{m_n} \, \delta_g^{j_n}(P) \, 
     P^{m_{n+1}} u\|_2
\leq c \, \|P^{m_1 + \ldots + m_{n+1} + n - j_1 - \ldots - j_n} u\|_2
\]
for all $u \in C^\infty(\Gamma)$ and 
$g \in W_{m_1 + \ldots + m_{n+1} + n}$.

In particular, if $m_1 + \ldots + m_{n+1} + n = j_1 + \ldots + j_n$,
then the operator 
\[
P^{m_1} \, \delta_g^{j_1}(P) \, P^{m_2} \ldots P^{m_n} \, \delta_g^{j_n}(P) \, 
     P^{m_{n+1}}
\]
extends to a bounded operator from $L_2(\Gamma)$ into $L_2(\Gamma)$.
\end{lemma}
\proof\
The proof is by induction to $n$.
The case $n = 1$ is done in Lemma~\ref{lpbdton603}.

If $m_1 + m_2 + 1 \geq j_1$ then it follows from Lemma~\ref{lpbdton603}
that 
\begin{eqnarray*}
\lefteqn{
\|P^{m_1} \, \delta_g^{j_1}(P) \, P^{m_2} \ldots P^{m_n} \, \delta_g^{j_n}(P) \, 
     P^{m_{n+1}} u\|_2
} \hspace*{30mm}  \\*
& \leq & c \, \|P^{m_1 + m_2 + 1 - j_1} \delta_g^{j_2}(P) \, P^{m_3} \, 
     \ldots P^{m_n} \, \delta_g^{j_n}(P) \, P^{m_{n+1}} u\|_2
\end{eqnarray*}
for a suitable constant $c$ and one can use the induction hypothesis.
Suppose that $m_1 + m_2 + 1 < j_1$.
Let $k \in \{ 2,\ldots,n \} $ be chosen minimal such that 
$m_1 + \ldots + m_{k+1} + k \geq j_1 + \ldots + j_k$.
Therefore $m_1 + \ldots + m_k + k - 1 < j_1 + \ldots + j_{k-1}$
and $m_{k+1} + 1 > j_k$.
Let $N = j_1 + \ldots + j_k - k - m_1 - \ldots - m_k$.
Then $N \in \{ 1,\ldots,m_{k+1} \} $.
Note that $m_1 + \ldots + m_k + k + N = j_1 + \ldots + j_k$.
Moreover, 
$N + m_k + 1 - j_k
= j_1 + \ldots + j_{k-1} - k + 1 - m_1 - \ldots - m_{k-1} > m_k \geq 0$.
So $N + m_k + 1 \geq j_k$.
Hence 
\begin{eqnarray*}
\lefteqn{
\|P^{m_1} \, \delta_g^{j_1}(P) \, P^{m_2} \ldots P^{m_n} \, \delta_g^{j_n}(P) \, 
     P^{m_{n+1}} u\|_2
} \hspace*{30mm}  \\*
& \leq & \|P^{m_1} \, \delta_g^{j_1}(P) \, P^{m_2} \ldots P^{m_{k-1}} \, 
     \delta_g^{j_{k-1}}(P) \, P^{m_k} \, \delta_g^{j_k}(P) \, P^N\|_{2 \to 2} \cdot  \\*
& & \hspace*{50mm} 
   {} \cdot 
   \|P^{m_{k+1} - N} \, \delta_g^{j_{k+1}}(P) \, \ldots \delta_g^{j_n}(P) \, P^{m_{n+1}} u\|_2  
.
\end{eqnarray*}
But by duality
\begin{eqnarray*}
\|P^{m_1} \ldots P^{m_{k-1}} \, 
     \delta_g^{j_{k-1}}(P) \, P^{m_k} \, \delta_g^{j_k}(P) \, P^N\|_{2 \to 2}  
& = & \|P^N \, \delta_g^{j_k}(P) \, P^{m_k} \, \delta_g^{j_{k-1}}(P) \, 
       \ldots \, P^{m_1}\|_{2 \to 2}  \\
& \leq & c \, \|P^{N + m_k - j_k + 1} \, \delta_g^{j_{k-1}}(P) \, 
      \ldots \, P^{m_1}\|_{2 \to 2} 
\end{eqnarray*}
for a suitable $c > 0$ by Lemma~\ref{lpbdton603}.
Now one can use twice the induction hypothesis.\hfill$\Box$

\begin{lemma} \label{lpbdton605}
Let $j,m \in \Ni$ and $k_1,k_2 \in \Ni_0$ with $k_1 + k_2 + m \geq j$.
Then there exists a $c > 0$ such that 
\[
\|P^{k_1} \, \delta_g^j(P^m) \, P^{k_2} u\|_2
\leq c \, \|P^{k_1 + k_2 + m - j} u\|_2
\]
for all $u \in C^\infty(\Gamma)$ and $g \in W_{k_1 + k_2 + m}$.
\end{lemma} 
\proof\
Since $\delta_g$ is a derivation, there are constants $c_{m_1,j_1,\ldots,m_{n+1}} \in \Ri$,
independent of $g$, such that 
\[
\delta_g^j(P^m)
= \sum c_{m_1,j_1,\ldots,m_{n+1}} \, P^{m_1} \, \delta_g^{j_1}(P) \, 
     P^{m_2} \, \delta_g^{j_2}(P) \ldots \delta_g^{j_n}(P) \, P^{m_{n+1}}
,  \]
where the sum is over all $n \in \{ 1,\ldots,j \} $,
$m_1,\ldots,m_{n+1} \in \Ni_0$ and $j_1,\ldots,j_n \in \Ni$ such that 
$j_1 + \ldots + j_n = j$ and $m_1 + \ldots + m_{n+1} + n = m$.
Now apply Lemma~\ref{lpbdton604}.\hfill$\Box$

\begin{lemma} \label{lpbdton606}
Let $k,m \in \Ni$, $\ell \in \Ni_0$ and $j_1,\ldots,j_k \in \{ 1,\ldots,m \} $.
Then there exists a $c > 0$ such that 
\[
\|P^\ell \, T^{(m)}_{z_{k+1}} \, \delta_g^{j_k}(P^m) \, T^{(m)}_{z_k} \, 
    \ldots T^{(m)}_{z_2} \delta_g^{j_1}(P^m) \, T^{(m)}_{z_1} u\|_2
\leq c \, \|P^{\ell + km - j_1 - \ldots - j_k} u\|_2
\]
for all $u \in C^\infty(\Gamma)$, $g \in W_{km + \ell}$
and $z_1,\ldots,z_{k+1} \in \Ci$
with $\RRe z_n > 0$ for all $n \in \{ 1,\ldots,k+1 \} $.
\end{lemma} 
\proof\
Since $T^{(m)}_z$ commutes with $P$ and $\|T^{(m)}_z\|_{2 \to 2} \leq 1$
this follows easily by 
induction from Lemma~\ref{lpbdton605}.\hfill$\Box$

\begin{lemma} \label{lpbdton607}
Let $m,n \in \Ni$, $\ell \in \Ni_0$ and $p \in [2,\infty]$.
Suppose that $d-1 < 2m$ and $n \leq m$.
Then there exists a $c > 0$ such that 
\[
\|P^\ell \, \delta_g^n(T^{(m)}_z)\|_{1 \to p}
\leq c \, (\RRe z)^{- \frac{d-1}{m} (1 - \frac{1}{p})} \,
       (\RRe z)^{-(\ell - n)/m} \, (\RRe z)^{-n} \, |z|^n 
\]
for all $z \in \Ci$ and $g \in W_{nm + \ell}$ with $\RRe z > 0$.
\end{lemma}
\proof\
We use Lemma~\ref{lpbdton402} to rewrite $\delta_g^n(T^{(m)}_z)$.
Let $k \in \{ 1,\ldots,n \} $ and $j_1,\ldots,j_k \in \Ni$ with 
$j_1 + \ldots + j_k = n$.
Let $(t_1,\ldots,t_{k+1}) \in H_k$.
There exists a $K \in \{ 1,\ldots,k+1 \} $ such that $t_K \geq \frac{1}{k+1}$.
Then 
\begin{eqnarray}
\lefteqn{
|z|^k \, \|P^\ell \, T^{(m)}_{t_{k+1} \, z} \, \delta_g^{j_k}(P^m) \, T^{(m)}_{t_k \, z} 
   \ldots
     T^{(m)}_{t_2 \, z} \, \delta_g^{j_1}(P^m) \, T^{(m)}_{t_1 \, z}\|_{1 \to p}
} \hspace*{10mm}  \nonumber \\*
& \leq & |z|^k \, \|T^{(m)}_{t_{k+1} \, z/2}\|_{2 \to p} 
  \cdot {} \nonumber  \\*
& & \hspace{1mm} {} \cdot
   \|P^\ell \, T^{(m)}_{t_{k+1} \, z/2} \, \delta_g^{j_k}(P^m) \, T^{(m)}_{t_k \, z} 
     \ldots 
     T^{(m)}_{t_2 \, z} \, \delta_g^{j_1}(P^m) \, T^{(m)}_{t_1 \, z/2}\|_{2 \to 2}
\|T^{(m)}_{t_1 \, z/2}\|_{1 \to 2}
\hspace*{5mm}
\label{elpbdton607;1}
\end{eqnarray}
By Lemma~\ref{lpbdton602} and duality there exists a suitable $c_1 > 0$ such that 
\begin{equation}
\|T^{(m)}_{t_{k+1} \, z/2}\|_{2 \to p}
\leq c_1 \, t_{k+1}^{- \frac{d-1}{m} (\frac{1}{2} - \frac{1}{p})} \, 
     (\RRe z)^{- \frac{d-1}{m} (\frac{1}{2} - \frac{1}{p})}
\quad \mbox{and} \quad
\|T^{(m)}_{t_1 \, z/2}\|_{1 \to 2}
\leq c_1 \, t_1^{-\frac{d-1}{2m}} \, (\RRe z)^{-\frac{d-1}{2m}}
.
\label{elpbdton607;2}
\end{equation}
We next estimate the big factor in (\ref{elpbdton607;1}).

Suppose that $K \in \{ 2,\ldots,k \} $.
Then 
\begin{eqnarray}
\lefteqn{
   \|P^\ell \, T^{(m)}_{t_{k+1} \, z/2} \, \delta_g^{j_k}(P^m) \, T^{(m)}_{t_k \, z} 
     \ldots 
     T^{(m)}_{t_2 \, z} \, \delta_g^{j_1}(P^m) \, T^{(m)}_{t_1 \, z/2}\|_{2 \to 2}
} \hspace*{25mm}  \nonumber \\*
& \leq & \|P^\ell \, T^{(m)}_{t_{k+1} \, z/2} \, \delta_g^{j_k}(P^m) \, 
       T^{(m)}_{t_k \, z} \ldots T^{(m)}_{t_{K+1} \, z} \,
                           \delta_g^{j_K}(P^m) \, T^{(m)}_{t_{K \, z/2}}\|_{2 \to 2}
       \cdot {}  \nonumber \\*
& & \hspace*{15mm} {} \cdot
   \|T^{(m)}_{t_1 \, \overline z/2} \, \delta_g^{j_1}(P^m) \, 
              T^{(m)}_{t_2 \, \overline z} \ldots T^{(m)}_{t_{K-1} \, \overline z} \,
                  \delta_g^{j_{K-1}}(P^m) \, T^{(m)}_{t_K \, \overline z/2}\|_{2 \to 2}
\hspace*{10mm}
\label{elpbdton607;3}
\end{eqnarray}
where we used duality in the second factor.
By Lemma~\ref{lpbdton606} and the decomposition
$T^{(m)}_{t_K \, z/2} = T^{(m)}_{t_K \, z/4} \circ T^{(m)}_{t_K \, z/4}$
there are suitable $c_2,c_3 > 0$ such that 
\begin{eqnarray*}
\lefteqn{
\|P^\ell \, T^{(m)}_{t_{k+1} \, z/2} \, \delta_g^{j_k}(P^m) \, 
       T^{(m)}_{t_k \, z} \ldots T^{(m)}_{t_{K+1} \, z} \,
                           \delta_g^{j_K}(P^m) \, T^{(m)}_{t_{K \, z/2}}\|_{2 \to 2}
} \hspace*{20mm}  \\*
& \leq & c_2 \, \|P^{\ell + (k-K+1)m - j_K - \ldots - j_k} \, T^{(m)}_{t_K \, z/4}\|_{2 \to 2}  \\
& \leq & c_3 \, (t_K \, \RRe z)^{-(\ell + (k-K+1)m - j_K - \ldots - j_k)/m} ,
\end{eqnarray*}
where we used the spectral theorem in the last step.
The second factor in (\ref{elpbdton607;3}) can be bounded similarly.
Since $t_K \geq \frac{1}{k+1}$, there is a suitable $c_4 > 0$ such that 
\[
\|P^\ell \, T^{(m)}_{t_{k+1} \, z/2} \, \delta_g^{j_k}(P^m) \, T^{(m)}_{t_k \, z} 
     \ldots 
     T^{(m)}_{t_2 \, z} \, \delta_g^{j_1}(P^m) \, T^{(m)}_{t_1 \, z/2}\|_{2 \to 2}
\leq c_4 \, (\RRe z)^{-(\ell + km - n)/m}
.  \]
Combining this with (\ref{elpbdton607;1}) and (\ref{elpbdton607;2})
one deduces that 
\begin{eqnarray*}
\lefteqn{
|z|^k \, \|P^\ell \, T^{(m)}_{t_{k+1} \, z} \, \delta_g^{j_k}(P^m) \, T^{(m)}_{t_k \, z} 
   \ldots
     T^{(m)}_{t_2 \, z} \, \delta_g^{j_1}(P^m) \, T^{(m)}_{t_1 \, z}\|_{1 \to p}
} \hspace*{40mm} \\*
& \leq & c_1^2 \, c_4 \, (\RRe z)^{- \frac{d-1}{m} (1 - \frac{1}{p})}
    (\RRe z)^{-(\ell + km - n)/m} \, |z|^k \, 
     t_1^{- \frac{d-1}{2m}} \, t_{k+1}^{-\frac{d-1}{m} (\frac{1}{2} - \frac{1}{p})}
. 
\end{eqnarray*}
The cases $K = 1$ and $K = k+1$ are similar.
Integrating over $H_k$ and taking the finite sum gives the result.\hfill$\Box$

\begin{lemma} \label{lpbdton608}
Let $m \in \Ni$ with $m \geq d$ and $\nu \in (0,1)$.
Then there exists a $c > 0$ such that 
\[
\|\delta_g^d(T^{(m)}_z)\|_{L_1 \to C^\nu}
\leq c \, (\cos \theta)^{(1 - \nu) / m} \,
        (\cos \theta)^{-d} \, |z|^{(1 - \nu) / m}
\]
for all $z \in \Ci$ and $g \in W_{dm + 1}$ with $\RRe z > 0$,
where $\theta = \arg z$.
\end{lemma}
\proof\
This follows from Lemma~\ref{lpbdton607} with $p = \frac{d-1}{1-\nu}$, $\ell = 1$ and 
$n = d$, followed by the Sobolev embedding of Proposition~\ref{ppbdtonA04}.\hfill$\Box$

\vertspace

At this stage we have the required bound for 
$\delta^d(T^{(m)}_z)$ from $L_1$ to $C^\nu$.
In order to obtain a bound for $\delta^d(S_z) = e^z \, \delta^d(T^{(1)}_z)$
we need a lemma on subordination.

\begin{lemma} \label{lpbdton609}
Let $-A$ be the generator of a semigroup in a Banach space $E$
which is bounded holomorphic in the sector $\Sigma_{\pi/2}^\circ$.
Let $F \in \cl(\cl(E))$ and $D \subset E$ a subspace.
Let $\cx,\cy$ be two Banach spaces with $D \subset \cx$.
Let $N \in \Ri$ and $\beta \in (-\infty,\frac{1}{2})$.
Suppose that $F(e^{-zA})u, F(e^{-z \sqrt{A}})u \in \cy$ and 
\[
\|F(e^{-zA}) u\|_\cy
\leq M \, (\cos \theta)^{-N} \, |z|^\beta \, \|u\|_\cx
\]
for all $u \in D$ and $z \in \Sigma_{\pi/2}^\circ$, where $\theta = \arg z$.
Then 
\[
\|F(e^{-z \sqrt{A}}) u\|_\cy
\leq c_\beta \, M \, (\cos \theta)^{-N} (\cos \theta)^{- (\frac{1}{2} - \beta)} \, 
   |z|^{2 \beta} \, \|u\|_\cx
\]
for all $u \in D$ and $z \in \Sigma_{\pi/2}^\circ$, 
where 
$c_\beta = \int_0^\infty \frac{1}{\sqrt{4 \pi}} \, s^{-3/2} \, e^{-\frac{1}{4s}} 
      \, s^\beta \, ds$ and $\theta = \arg z$.
\end{lemma}
\proof\
For all $z \in \Ci$ with $\RRe z > 0$ define $\mu_z \colon (0,\infty) \to \Ci$
by 
\[
\mu_z(s) = \frac{1}{\sqrt{4 \pi}} \, z \, s^{-3/2} \, e^{-\frac{z^2}{4s}}
.  \]
Then 
\[
e^{-t \sqrt{B}} = \int_0^\infty \mu_t(s) \, e^{-s B} \, ds
\]
for all $t > 0$ and every bounded strongly continuous semigroup 
by the example on page 268 in \cite{Yos}.
Fix $z = r \, e^{i \theta} \in \Ci$ with $|\theta| < \frac{\pi}{2}$ and $r \in (0,\infty)$.
Choosing $B = e^{i \theta} A$ gives
\begin{equation}
e^{- t e^{i \theta/2} \sqrt{A}}
= e^{-t \sqrt{B}}
= \int_0^\infty \mu_t(s) \, e^{-s e^{i \theta} A} \, ds
\label{elpbdton609;1}
\end{equation}
for all $t \in (0,\infty)$.
Since both sides in (\ref{elpbdton609;1}) extend holomorphically 
to the sector $\Sigma_{\pi/4}^\circ$
one deduces that 
\[
e^{-z \sqrt{A}}
= \int_0^\infty \mu_{r e^{i \theta/2}}(s) \, e^{-s e^{i \theta} A} \, ds
.  \]
Now let $u \in D$.
Then 
\begin{eqnarray*}
\|F(e^{-z \sqrt{A}}) u\|_\cy
& \leq & \int_0^\infty |\mu_{r e^{i \theta/2}}(s)| \, \|F(e^{-s e^{i \theta} A}) u\|_\cy \, ds  \\
& \leq & \int_0^\infty |\mu_{r e^{i \theta/2}}(s)| \, 
     M \, (\cos \theta)^{-N} \, s^\beta \, \|u\|_\cx \, ds
.  
\end{eqnarray*}
But 
\begin{eqnarray*}
\int_0^\infty |\mu_{r e^{i \theta/2}}(s)| \, s^\beta \, ds
& \leq & \int_0^\infty \frac{1}{\sqrt{4 \pi}} \, r \, s^{-3/2} \, e^{-\frac{r^2 \cos \theta}{4s}} 
      \, s^\beta \, ds  \\
& = & r^{2 \beta} (\cos \theta)^{\beta - \frac{1}{2}} 
    \int_0^\infty \frac{1}{\sqrt{4 \pi}} \, s^{-3/2} \, e^{-\frac{1}{4s}} 
      \, s^\beta \, ds  \\
& = & c_\beta \, r^{2 \beta} (\cos \theta)^{\beta - \frac{1}{2}} 
.
\end{eqnarray*}
Therefore 
\[
\|F(e^{-z \sqrt{A}}) u\|_\cy
\leq c_\beta \, M \, (\cos \theta)^{-N + \beta - \frac{1}{2}} \, r^{2 \beta} \, \|u\|_\cx
\]
as required.\hfill$\Box$

\begin{lemma} \label{lpbdton610}
Let $\nu \in (0,1)$.
Then there exists a $c > 0$ such that 
\[
\|\delta_g^d(T^{(1)}_z)\|_{L_1 \to C^\nu}
\leq c \, (\cos \theta)^{1 - \nu} \, (\cos \theta)^{-d-k/2} \, |z|^{1 - \nu}
\]
for all $z \in \Ci$ and $g \in W_{2^k d + 1}$ with $\RRe z > 0$,
where $k = \lceil \frac{\log d}{\log 2} \rceil$ and $\theta = \arg z$.
\end{lemma}
\proof\
Note that $2^k \geq d$.
Let $c > 0$ be as in Lemma~\ref{lpbdton608} with $m = 2^k$.
For all $\beta \in (0,\frac{1}{2})$ let $c_\beta$ be as in Lemma~\ref{lpbdton609}.
Using Lemma~\ref{lpbdton609} it follows by induction to $\ell$ that 
\begin{eqnarray*}
\lefteqn{
\|\delta_g^d(T^{(2^{k-\ell})}_z)\|_{L_1 \to C^\nu}
} \hspace*{10mm}  \\*
& \leq & c \, c_{(1 - \nu) / 2^k} \ldots c_{(1 - \nu)/2^{k-\ell+1}} \,
   (\cos \theta)^{-d + \frac{1-\nu}{2^k}} \, 
   (\cos \theta)^{-(\frac{1}{2} - \frac{1-\nu}{2^k})} \ldots 
        (\cos \theta)^{-(\frac{1}{2} - \frac{1-\nu}{2^{k - \ell + 1}})} \, 
   |z|^{\frac{1-\nu}{2^{k - \ell}}}  \\
& = & c \, c_{(1 - \nu) / 2^k} \ldots c_{(1 - \nu)/2^{k-\ell+1}} \, 
     (\cos \theta)^{-d - \ell / 2} \, (\cos \theta)^{\frac{1-\nu}{2^{k-\ell}}} \,
       |z|^{\frac{1-\nu}{2^{k - \ell}}} 
\end{eqnarray*}
for all $\ell \in \{ 0,\ldots,k \} $.
Choosing $\ell = k$ gives the estimate of the lemma.\hfill$\Box$

\vertspace

We are now able to prove the main theorem of this section.

\begin{thm} \label{tpbdton611}
For all $\nu \in (0,1)$ there exists a $c > 0$ such that 
\[
|K_z(x,y)|
\leq c \, (\cos \theta)^{1-\nu} \, (\cos \theta)^{-d-k/2} \, 
    \frac{(|z| \wedge 1)^{-(d-1)}}{\displaystyle \Big( 1 + \frac{|x-y|}{|z|} \Big)^{d-\nu} }
   \, 
\]
for all $z \in \Ci$ with $\RRe z > 0$,
where $k = \lceil \frac{\log d}{\log 2} \rceil$ and $\theta = \arg z$.
\end{thm}
\proof\
Let $c > 0$ be as in Lemma~\ref{lpbdton610}.
By Lemma~\ref{lpbdton300.6} there exists a $c_0 > 0$ such that 
\[
\frac{1}{c_0} \, \rho_\Gamma(x,y)
\leq \rho^{(2^k d)}_\Gamma(x,y)
\leq c_0 \, \rho_\Gamma(x,y)
\]
for all $x,y \in \Gamma$.
Let $g \in W_{2^k d + 1}$.
Then 
\[
|(\delta_g^d(T^{(1)}_z) u)(x) - (\delta_g^d(T^{(1)}_z) u)(x')|
\leq c \, (\cos \theta)^{1-\nu} \, (\cos \theta)^{-d-k/2} \, |z|^{1-\nu} \, \rho_\Gamma(x,x')^\nu \, \|u\|_1
\]
for all $u \in L_1(\Gamma)$ and $x,x' \in \Gamma$.
Hence 
\begin{eqnarray*}
\lefteqn{
|(g(x) - g(y))^d \, K_z(x,y) 
     - (g(x') - g(y))^d \, K_z(x',y)| \, e^{-\RRe z}
} \hspace{50mm} \\*
& \leq & c \, (\cos \theta)^{1-\nu} \, (\cos \theta)^{-d-k/2} \, |z|^{1-\nu} \, \rho_\Gamma(x,x')^\nu
\end{eqnarray*}
for all $x,x',y \in \Gamma$.
Choosing $x' = y$ gives 
\[
|g(x) - g(y)|^d \, |K_z(x,y)| \, e^{-\RRe z}
\leq c \, (\cos \theta)^{1-\nu} \, (\cos \theta)^{-d-k/2} \, |z|^{1-\nu} \, \rho_\Gamma(x,x')^\nu
\]
for all $x,y \in \Gamma$.
Optimizing over $g \in W_{2^k d}$ it follows that
\[
c_0^{-d} \, \rho_\Gamma(x,y)^d \, |K_z(x,y)| \, e^{-\RRe z}
\leq c \, (\cos \theta)^{1-\nu} \, (\cos \theta)^{-d-k/2} \, |z|^{1-\nu} \, \rho_\Gamma(x,x')^\nu
\]
for all $x,y \in \Gamma$.
Therefore 
\begin{equation}
\Big( \frac{\rho_\Gamma(x,y)}{|z|} \Big)^{d-\nu} |K_z(x,y)|
\leq c \, c_0^d \, (\cos \theta)^{1-\nu} \, (\cos \theta)^{-d-k/2} \, |z|^{-(d-1)} \, e^{\RRe z}
\label{etpbdton611;1}
\end{equation}
for all $x,y \in \Gamma$.
It follows from Lemma~\ref{lpbdton602} and duality that there exists a 
suitable $c_1 > 0$ such that 
\begin{equation}
|K_z(x,y)| 
\leq \|T^{(1)}_z\|_{1 \to \infty} \, e^{\RRe z}
\leq c_1 (\cos \theta)^{-(d-1)} \, |z|^{-(d-1)} \, e^{\RRe z}
\label{etpbdton611;2}
\end{equation}
for all $x,y \in \Gamma$.
Then the theorem for $|z| \leq 1$ follows 
from adding (\ref{etpbdton611;1}) and (\ref{etpbdton611;2}), 
together with Lemma~\ref{lpbdton300.5}\ref{lpbdton300.5-4}.

Finally we deal with the case $|z| \geq 1$.  
Let $C > 0$ be as in Theorem \ref{th2.5}.
Then for all $z = t + is$ with $ t > 0$ one estimates
\begin{eqnarray}
\|S_z\|_{1 \to \infty}
& \leq & \|S_{t/2}\|_{2 \to \infty} \, \|S_{is}\|_{2 \to 2} \, \|S_{t/2}\|_{1 \to 2} \nonumber \\
& \leq & C^2  \Big(1 \wedge  \tfrac{1}{2} |z| \cos \theta \Big)^{-(d-1) }
\leq 2^{d-1} C^2 \, (\cos \theta)^{-(d-1)}  (|z| \wedge 1 )^{-(d-1)}. 
\label{etpbdton611;3}
\end{eqnarray}
Since $\Gamma$ is bounded, there exists a $c > 0$ such that 
\[
|K_z(x,y)|
\leq c \, \frac{\|S_z\|_{1 \to \infty}}{\displaystyle \Big( 1 + \frac{|x-y|}{|z|} \Big)^{d-\nu} }
\leq 2^{d-1} c \, C^2 \, (\cos \theta)^{-(d-1)} \, 
    \frac{(|z| \wedge 1 )^{-(d-1)}}{\displaystyle \Big( 1 + \frac{|x-y|}{|z|} \Big)^{d-\nu} }
\]
for all $x,y \in \Gamma$ and $z \in \Ci$ with $\RRe z > 0$ and $|z| \geq 1$.
This completes the proof of the theorem.\hfill$\Box$

\begin{cor} \label{cpbdton612}
For all $\nu \in (0,1)$ there exists a $c > 0$ such that 
\[
\|S_z\|_{p \to p}
\leq c \,  (\cos \theta)^{-d-k/2 + 1-\nu} \, 
\]
for all $p \in [1,\infty]$ and $z \in \Ci$ with $\RRe z > 0$,
where $k = \lceil \frac{\log d}{\log 2} \rceil$ and $\theta = \arg z$.
\end{cor}
\proof\
The bounds for $p = 1$ follows from a quadrature estimate from the 
Poisson bounds in Theorem~\ref{tpbdton611}.
Then the bounds for  $p \in (1,\infty]$ follow from duality and 
interpolation.\hfill$\Box$

\begin{cor} \label{cpbdton613}
For all $\nu \in (0,1)$ there exists a $c > 0$ such that 
\[
\|S_z\|_{p \to q}
\leq c \, (\cos \theta)^{-d-k/2 + 1-\nu} \, |z|^{-(d-1)(\frac{1}{p} - \frac{1}{q})}
\]
for all $p,q \in [1,\infty]$ and $z \in \Ci$ with $p \leq q$,
$\RRe z > 0$ and $|z| \leq 1$,
where $k = \lceil \frac{\log d}{\log 2} \rceil$ and $\theta = \arg z$.
\end{cor}
\proof\
This follows from interpolation of the bounds of Corollary~\ref{cpbdton612}
and the bounds (\ref{etpbdton611;3}).\hfill$\Box$

\vertspace

We are finally able to prove the full Poisson bounds for complex $z$.

\vertspace

\noindent
{\bf Proof of Theorem~\ref{tpbdton102}\hspace*{5pt}\ }
This follows from Proposition~\ref{ppbdton401} and Corollary~\ref{cpbdton613},
similarly as in the proof
of Theorem~\ref{tpbdton611}.\hfill$\Box$

\section{Derivatives} \label{Spbdton7}

The kernel $K_z$ of the operator $S_z$ is a smooth function. 
The aim of this section is to prove Poisson bounds for the 
spacial derivatives of $K_z$. 
If confusion is possible, then we denote by a subscript 
$(1)$ and $(2)$ the first or second variable on which 
a derivative acts.

The main theorem of this section is the following.

\begin{thm} \label{tpbdton701}
For all $k,\ell \in \Ni_0$ there exists a $c > 0$ such that 
\[
|(\nabla_{(1)}^k \nabla_{(2)}^\ell K_z)(x,y)|
\leq c \, (\cos \theta)^{-4d(d+1) - k - \ell} \, 
    \frac{|z|^{-(d-1)} \, |z|^{-(k + \ell)} \, e^{2 |z|} } 
         {\displaystyle \Big( 1 + \frac{|x-y|}{|z|} \Big)^d }  
\]
for all $z \in \Ci$ and $x,y \in \Gamma$ with $\RRe z > 0$,
where $\theta = \arg z$.
\end{thm}

The proof uses interpolation and the Poisson bounds of Theorem~\ref{tpbdton102}.
The first step is that Theorem~\ref{tpbdton102} has an easy corollary.

\begin{cor} \label{cpbdton702}
There exists a $c_0 > 0$ such that 
\[
\|\delta_g^j(S_z)\|_{1 \to \infty}
\leq c_0 \, (\cos \theta)^{-2d(d+1)} \, 
    (|z| \wedge 1)^{-(d-1)} \, |z|^j 
\]
for all $j \in \{ 0,\ldots,d \} $ and $z \in \Ci$ with $\RRe z > 0$,
where $\theta = \arg z$.
\end{cor}

The key estimate for the proof of Theorem~\ref{tpbdton701} is in the 
following lemma.

\begin{lemma} \label{lpbdton703}
Let $(V,\varphi)$ be a chart, $\chi \in C_c^\infty(V)$, $j \in \{ 0,\ldots,d \} $
and $\alpha$ a multi-index over $ \{ 1,\ldots,d-1 \} $.
Then there exists a $c > 0$ such that 
\[
\|\Big( \frac{\partial}{\partial \varphi} \Big)^\alpha M_\chi \, 
      \delta_g^j(T^{(1)}_z)\|_{1 \to \infty}
\leq c \, (|z| \wedge 1)^{-(d-1)} \, |z|^j \, |z|^{-|\alpha|} \,
   (\cos \theta)^{-2d(d+2) - |\alpha|}
\]
for all $g \in W_{2^k j + d + |\alpha| + 1}$
and $z \in \Ci$ with $\RRe z > 0$,
where $\theta = \arg z$.
\end{lemma}
\proof\
Let $k = \lceil \frac{\log d}{\log 2} \rceil$ and $\ell = d + |\alpha| + 1$. 
Let $p \in (d-1,\infty)$.
By Lemma~\ref{lpbdton607} there exists a $c_1 > 0$ such that 
\[
\|P^\ell \, \delta_g^j(T^{(2^k)}_z)\|_{1 \to p}
\leq c_1 \, |z|^{-\frac{d-1}{2^k}(1 - \frac{1}{p})} \, |z|^{-\frac{(\ell-j)}{2^k}}
   \, (\cos \theta)^{-2d - \ell}
\]
for all $g \in W_{2^k j + \ell}$ and $z \in \Ci$ with $\RRe z > 0$.
Arguing as in the proof of Lemma~\ref{lpbdton610} one deduces that 
there exists a $c_2 > 0$ such that 
\[
\|P^\ell \, \delta_g^j(T^{(1)}_z)\|_{1 \to p}
\leq c_2 \, |z|^{-(d-1)(1 - \frac{1}{p})} \, |z|^{-(\ell - j)}
   \, (\cos \theta)^{-2d - \ell - k}
\]
for all $g \in W_{2^k j + \ell}$ and $z \in \Ci$ with $\RRe z > 0$.
Hence by Proposition~\ref{ppbdton601} there exists a $c_3 > 0$ such that 
\[
\|\delta_g^j(T^{(1)}_z)\|_{L_1(\Gamma) \to W^{\ell,p}(\Gamma)}
\leq c_3 \, |z|^{-(d-1)(1 - \frac{1}{p})} \, |z|^{-(\ell - j)}
   \, (\cos \theta)^{-2d - \ell - k}
\]
for all $g \in W_{2^k j + \ell}$ and $z \in \Ci$ with $\RRe z > 0$.
Next, let $c_0 > 0$ be as in Corollary~\ref{cpbdton702}.
By Proposition~\ref{ppbdtonA02.2} there exists a $c_4 > 0$ such that 
\[
\|\Big( \frac{\partial}{\partial \varphi} \Big)^\alpha (\chi \, u)\|_{L_\infty(\Gamma)}
\leq c_4 \, \|\chi \, u\|_{W^{\ell,p}(\Gamma)}^\gamma \, 
            \|\chi \, u\|_{L_\infty(\Gamma)}^{1 - \gamma}
\]
for all $u \in C^\infty(\Gamma)$, where 
\[
\gamma = \frac{|\alpha|}{\ell - \frac{d-1}{p}}
.  \]
Therefore 
\begin{eqnarray*}
\lefteqn{
\|\Big( \frac{\partial}{\partial \varphi} \Big)^\alpha M_\chi \, 
           \delta_g^j(T^{(1)}_z)\|_{1 \to \infty}
} \hspace{10mm} \\*
& \leq & c_4 \, \|M_\chi \, \delta_g^j(T^{(1)}_z)\|_{L_1(\Gamma) \to W^{\ell,p}(\Gamma)}^\gamma \, 
            \|M_\chi \, \delta_g^j(T^{(1)}_z)\|_{L_1(\Gamma) \to L_\infty(\Gamma)}^{1 - \gamma}  \\
& \leq & c_4 \, \Big( c_3 \, \|M_\chi\|_{W^{\ell,p}(\Gamma) \to W^{\ell,p}(\Gamma)} \,
     |z|^{-(d-1)(1 - \frac{1}{p})} \, |z|^{-(\ell - j)} \,
    (\cos \theta)^{2d - \ell - k}
                \Big)^\gamma  \cdot  \\*
& & \hspace{50mm} {} \cdot
   \Big( c_0 \, \|\chi\|_\infty \, (|z| \wedge 1)^{-(d-1)} \, |z|^j \, (\cos \theta)^{-2d(d+1)} 
   \Big)^{1 - \gamma}  \\
& \leq & c_5 \, (|z| \wedge 1)^{-(d-1)} \, |z|^{-|\alpha|} \, |z|^j \, (\cos \theta)^{-2d(d+1) - \ell}
\end{eqnarray*}
for a suitable constant $c_5$.\hfill$\Box$

\begin{lemma} \label{lpbdton704}
Let $(V,\varphi)$ be a chart, $\chi \in C_c^\infty(V)$,
and $\alpha$ a multi-index over $ \{ 1,\ldots,d-1 \} $.
Then there exists a $c > 0$ such that 
\[
\|\delta_g^d( \Big( \frac{\partial}{\partial \varphi} \Big)^\alpha M_\chi \, T^{(1)}_z)\|_{1 \to \infty}
\leq c \, (|z| \wedge 1)^{1-|\alpha|} \, (|z| \vee 1)^d \, 
    (\cos \theta)^{-2d(d+2) - |\alpha|} \, 
\]
for all $g \in W_{2^k j + d + |\alpha| + 1}$
and $z \in \Ci$ with $\RRe z > 0$,
where $\theta = \arg z$.
\end{lemma}
\proof\
It follows by induction to $m$ that for all $m \in \Ni_0$ and 
multi-indices $\beta_1,\ldots,\beta_m,\gamma$ over $ \{ 1,\ldots,d-1 \} $
there is a constant $c_{\beta_1,\ldots,\beta_m,\gamma} \in \Ri$ such that 
\begin{equation}
\delta_g^m( \Big( \frac{\partial}{\partial \varphi} \Big)^\alpha M_\chi )
= \sum c_{\beta_1,\ldots,\beta_m,\gamma} \, 
   M_{ (\frac{\partial}{\partial \varphi} )^{\beta_1} g } 
     \ldots 
   M_{ (\frac{\partial}{\partial \varphi} )^{\beta_m} g } 
    \Big( \frac{\partial}{\partial \varphi} \Big)^\gamma \, M_\chi
\label{elpbdton704;1}
\end{equation}
uniformly for all $g \in C^\infty(\Gamma)$, where the sum 
is over all $\beta_1,\ldots,\beta_m,\gamma$ with 
$|\beta_1|, \ldots, |\beta_m| \in \Ni$ and 
$|\beta_1| + \ldots + |\beta_m| + |\gamma| = |\alpha|$.
Note that $|\alpha| - |\gamma| \geq m$.
Since $\delta_g$ is a derivation, one has 
\[
\delta_g^d( \Big( \frac{\partial}{\partial \varphi} \Big)^\alpha M_\chi \, T^{(1)}_z)
= \sum_{j=0}^d {d \choose j} \,
   \delta_g^{d-j}(\Big( \frac{\partial}{\partial \varphi} \Big)^\alpha M_\chi ) \,
   \delta_g^j(T^{(1)}_z)
.  \]
Now use (\ref{elpbdton704;1}) and Lemma~\ref{lpbdton703}.\hfill$\Box$

\begin{lemma} \label{lpbdton705}
Let $(V,\varphi)$ be a chart, $\chi \in C_c^\infty(V)$,
and $\alpha$ a multi-index over $ \{ 1,\ldots,d-1 \} $.
Then there exists a $c > 0$ such that 
\[
|\Big( 
\Big( \frac{\partial}{\partial \varphi} \Big)^\alpha_{(1)}
   ( (\chi \otimes \one) K_z ) 
\Big) (x,y)|
\leq c \, \frac{ |z|^{-(d-1)} \, |z|^{-|\alpha|} \, e^{2|z|} } 
         {\displaystyle \Big( 1 + \frac{|x-y|}{|z|} \Big)^d }
   \, (\cos \theta)^{-2d(d+2) - |\alpha|}
\]
for all $t > 0$ and $x,y \in M$.
\end{lemma}
\proof\
This follows from Lemma~\ref{lpbdton704} by minimizing over $g$,
together with the bounds of Lemma~\ref{lpbdton703} with $j = 0$.\hfill$\Box$

\vertspace

In order to have derivatives on both variables we use duality and the 
next lemma, which states that the convolution of two Poisson bounds 
is again a Poisson bound.

\begin{lemma} \label{lpbdton706}
There exists a $c > 0$ such that 
\[
\int_\Gamma \frac{ (t \wedge 1)^{-(d-1)} }
            {\displaystyle \Big( 1 + \frac{|x-z|}{t} \Big)^d }
  \cdot \frac{ (t \wedge 1)^{-(d-1)} }
             {\displaystyle \Big( 1 + \frac{|z-y|}{t} \Big)^d }
  \, dz
\leq c \, \frac{ (t \wedge 1)^{-(d-1)} }
               {\displaystyle \Big( 1 + \frac{|x-y|}{t} \Big)^d }
\]
for all $t > 0$ and $x,y \in \Gamma$.
\end{lemma}
\proof\
For all $t > 0$ define $P_t \colon \Gamma \times \Gamma \to \Ri$ by 
\[
P_t(x,y) 
= \int_\Gamma \frac{ (t \wedge 1)^{-(d-1)} }
              {\displaystyle \Big( 1 + \frac{|x-z|}{t} \Big)^d }
  \cdot \frac{ (t \wedge 1)^{-(d-1)} }
             {\displaystyle \Big( 1 + \frac{|z-y|}{t} \Big)^d }
  \, dz
.  \]
Let
\[
c_0 = \sup_{t \in (0,\infty)} \, 
    \sup_{x \in \Gamma} 
       \int_\Gamma \frac{ (t \wedge 1)^{-(d-1)} }
                   { \displaystyle \Big( 1 + \frac{|x-y|}{t} \Big)^d }
< \infty
.  \]
Let $t > 0$ and $x,y \in \Gamma$.
Then
\[
P_t(x,y) 
\leq \int_\Gamma \frac{ (t \wedge 1)^{-(d-1)} }
                 {\displaystyle \Big( 1 + \frac{|x-z|}{t} \Big)^d }
    \cdot (t \wedge 1)^{-(d-1)}
        \, dz
\leq c_0 \, (t \wedge 1)^{-(d-1)}
.  \]
Moreover, $|x-y|^d \leq (|x-z| + |z-y|)^d \leq 2^d (|x-z|^d + |z-y|^d)$.
Hence 
\[
|x-y|^d \, P_t(x,y)
\leq 2^d \int_\Gamma \frac{ (t \wedge 1)^{-(d-1)} }
              {\displaystyle \Big( 1 + \frac{|x-z|}{t} \Big)^d }
  \cdot (|x-z|^d + |z-y|^d) \cdot 
        \frac{ (t \wedge 1)^{-(d-1)} }
             {\displaystyle \Big( 1 + \frac{|z-y|}{t} \Big)^d }
  \, dz
.  \]
But 
\begin{eqnarray*}
\int_\Gamma \frac{ (t \wedge 1)^{-(d-1)} |x-z|^d}
              {\displaystyle \Big( 1 + \frac{|x-z|}{t} \Big)^d }
  \cdot 
        \frac{ (t \wedge 1)^{-(d-1)} }
             {\displaystyle \Big( 1 + \frac{|z-y|}{t} \Big)^d }
  \, dz
& \leq & \int_\Gamma  (t \wedge 1)^{-(d-1)} \, t^d \, 
        \frac{ (t \wedge 1)^{-(d-1)} }
             {\displaystyle \Big( 1 + \frac{|z-y|}{t} \Big)^d } 
          \, dz \\
& \leq & c_0 \, (t \wedge 1)^{-(d-1)} \, t^d 
.
\end{eqnarray*}
Estimating similarly the other term one deduces that 
\[
|x-y|^d \, P_t(x,y) \leq 2^{d+1} c_0 \, (t \wedge 1)^{-(d-1)} \, t^d
.  \]
Then the lemma follows with $c = (1 + 2^{d+1}) c_0$.\hfill$\Box$

\vertspace

\noindent
{\bf Proof of Theorem~\ref{tpbdton701}\hspace{5pt}}\
Let $(V_1,\varphi_1)$ and $(V_2,\varphi_2)$ be charts, $\chi_1 \in C_c^\infty(V_1,\Ri)$,
$\chi_2 \in C_c^\infty(V_2,\Ri)$,
and $\alpha$ and $\beta$ be multi-indices over $ \{ 1,\ldots,d-1 \} $.
The semigroup property gives
\begin{eqnarray*}
\lefteqn{
\Big( 
\Big( \frac{\partial}{\partial \varphi_1} \Big)^\alpha_{(1)}
\Big( \frac{\partial}{\partial \varphi_2} \Big)^\beta_{(2)}
   ( (\chi_1 \otimes \chi_2) K_{2z} ) 
\Big) (x,y)
} \hspace{20mm} \\*
& = & \int_\Gamma
     \Big( \Big( \frac{\partial}{\partial \varphi_1} \Big)^\alpha_{(1)}
   ( (\chi_1 \otimes \one) K_z ) 
     \Big) (x,x')
\cdot
     \Big( \Big( \frac{\partial}{\partial \varphi_2} \Big)^\alpha_{(2)}
   ( (\one \otimes \chi_2) K_z ) 
     \Big) (x',y)
   \, dx'
\end{eqnarray*}
for all $x,y \in \Gamma$ and $z \in \Ci$ with $\RRe z > 0$.
But 
\[
    \Big( \Big( \frac{\partial}{\partial \varphi_2} \Big)^\alpha_{(2)}
   ( (\one \otimes \chi_2) K_z ) 
     \Big) (x',y)
=     \Big( \Big( \frac{\partial}{\partial \varphi_2} \Big)^\alpha_{(1)}
   ( (\chi_2 \otimes \one) \overline{K_z} ) 
     \Big) (y,x')
.  \]
Using Lemmas~\ref{lpbdton705} and \ref{lpbdton706} it follow that there 
exists a $c > 0$ such that 
\[
|\Big( 
\Big( \frac{\partial}{\partial \varphi_1} \Big)^\alpha_{(1)}
\Big( \frac{\partial}{\partial \varphi_2} \Big)^\beta_{(2)}
   ( (\chi_1 \otimes \chi_2) K_z ) 
\Big) (x,y)|
\leq c \, \frac{ |z|^{-(d-1)} \, |z|^{-(|\alpha| + |\beta|)} \, e^{2|z|} } 
         {\displaystyle \Big( 1 + \frac{|x-y|}{|z|} \Big)^d }
     \, (\cos \theta)^{-4d(d+2) - |\alpha| - |\beta|}
\]
for all $x,y \in \Gamma$ and $z \in \Ci$ with $\RRe z > 0$.
Now the theorem follows by a partition of the unity and 
Lemma~\ref{lpbdtonA07}.\hfill$\Box$

\section{Holomorphy and $H_\infty$-functional calculus} \label{Spbdton8}

In this  section we give applications of our Poisson bounds to the $L_p$-holomorphy 
of the semigroup as well as $H_\infty$-functional calculus and sharp spectral multipliers. 
We start with the holomorphy.
Recall that the operator $\cn_V$
is self-adjoint and hence the semigroup $S^V$  is  holomorphic on the sector 
$\Sigma_{\pi/2}^\circ$ in $L_2(\Gamma)$, where $\Sigma_\alpha$ is defined in~(\ref{eSpbdton4;1}).
If $V \geq 0$  then $S^V_z$ is a contraction operator on $L_2(\Gamma)$ 
for every $z \in \Sigma_{\pi/2}^\circ$.
On the other hand, the Poisson bound we proved 
allow to extend the semigroup
$S^V$ from $L_p(\Gamma) \cap L_2(\Gamma)$ to  a strongly continuous semigroup on 
$L_p(\Gamma)$ for all $p \in [1, \infty)$.
A natural question concerns the holomorphy 
of the extension to $L_1(\Gamma$) and  describe the sector of holomorphy.
It is now well known (cf.\ \cite{Ouh4} or \cite{Ouh5} Corollary~7.5) that a Gaussian 
upper bound of the heat kernel of a self-adjoint  semigroup implies analyticity on 
$L_1$ on the sector $\Sigma_{\pi/2}^\circ$.
This fact  is not clear  if instead we have Poisson bounds.
Nevertheless we have the following result.

\begin{thm}\label{holomo}
Suppose $0 \leq V \in L_\infty(\Omega)$.
The semigroup $S^V$ is holomorphic on $L_1(\Gamma)$ on the sector 
$\Sigma_{\frac{\pi}{2d}}^\circ$.
If $V = 0$, then $S$ is holomorphic on $L_1(\Gamma)$ 
on the sector $\Sigma_{\pi/2}^\circ$.
\end{thm}
\proof\
For all $z \in \Ci$ with $\RRe z > 0$ let $K^V_z$ be the kernel of $S^V_z$.
By Theorem \ref{tpbdton101} and Proposition 3.3 in \cite{DR1} it follows 
that for all $\varepsilon \in (0,1)$ and $\theta \in (0, \varepsilon \frac{\pi}{2})$ 
there is a $C > 0$ such that 
\begin{equation}\label{ehol1}
| K^V_z (x,y) | 
\leq C  \frac{(1 \wedge \RRe z)^{-(d-1)}}{\displaystyle \Big( 1 + \frac{|x-y|}{|z|} \Big)^{d(1-\varepsilon)} }
\end{equation}
for all $z \in \Sigma_\theta^\circ$ and $x,y \in \Gamma$.
Now suppose hat $d \, \varepsilon < 1$.
Then by \cite{DR1} Proposition~2.3 the semigroup $t \mapsto S^V_{t \, e^{i \varphi}}$
extends to a $C_0$-semigroup on $L_1(\Gamma)$ for each $\varphi \in (-\theta,\theta)$.
Integrating the bounds of (\ref{ehol1}) on the $(d-1)$-dimensional manifold  $\Gamma$ we see that 
there is a $C' > 0$ such that
\[
\int_\Gamma | K^V_z (x,y) |  \, d\sigma(x) \le C'
\]
for all $z \in \Sigma_\theta^\circ$ and $y \in \Gamma$.
Therefore the semigroups $(S^V_{t \, e^{i \varphi}})_{t > 0}$ are bounded, uniformly 
for all $\varphi \in (-\theta,\theta)$.
Hence $S^V$ is holomorphic on $L_1(\Gamma)$ on the sector 
$\Sigma_\theta^\circ$ by \cite{Kat1} Theorem~IX.1.23.
This means that we have holomorphy of $S^V$ on $L_1(\Gamma)$ on the sector 
$\Sigma_{\frac{\pi}{2d}}^\circ$.

If $V = 0$ we apply Theorem \ref{tpbdton102} to obtain the second assertion.\hfill$\Box$

\vertspace

We do not know whether $S^V$ is holomorphic on the right half-plane on $L_1(\Gamma)$.
Another application of Theorem \ref{tpbdton101} concerns the $H_\infty$-functional calculus.

\begin{thm} \label{tpbdton802}
Suppose $V \geq 0$.   
Let $\mu \in (\frac{\pi(d-1)}{2d},\pi)$ and $p \in (1,\infty)$.
Then $\cn_V$ has a bounded $H_\infty(\Sigma_\mu^\circ)$-functional calculus 
on $L_p(\Gamma)$.
Moreover, $f(\cn_V)$ is of weak type $(1,1)$ for all $f \in H_\infty(\Sigma_\mu^\circ)$.

If $V = 0$ then the above is valid for all $\mu \in (0,\pi)$.
\end{thm}
\proof\
This follows from (\ref{ehol1}) and Theorem 3.1 in \cite{DR1} if $V \not= 0$.
If $V = 0$ we can use the bounds  for complex time in Theorem \ref{tpbdton102},
which allow any choice of $\mu \in (0,\pi)$.\hfill$\Box$

\vertspace

An interesting particular case of  the  holomorphic functional calculus is the 
boundedness on $L_p(\Gamma)$ of imaginary powers  $\cn_V^{is}$.
The bounded imaginary powers on $L_p(\Gamma)$ in case $V = 0$ were proved before by 
Escher--Seiler \cite{EscS} with different methods.

We emphasize that for the operator $\cn$,  stronger results are known.
Indeed  a spectral multiplier theorem is proved in \cite{SeeSog}, Theorem 3.1. 
More precisely, it  follows from the results there  that  $f(\cn)$ is bounded on 
$L_p(\Gamma)$ for all $p \in (1,\infty)$ provided $f$ satisfies the H\"ormander condition
\[
 \sup_{t > 0} \| f(.) \beta(t.) \|_{W^{2,s}} < \infty,
\]
where $\beta$ is a smooth non-trivial auxiliary function and $s > \frac{d-1}{2}$.
It follows easily from the Cauchy formula that the latter condition holds 
if $f$ is a bounded holomorphic function in some sector of angle $\mu > 0$. 

Note that using our Poisson bound one can adapt the method from \cite{DOS}  
to obtain the previously mentioned spectral multiplier result for $\cn$.
Indeed, if one uses 
Theorem \ref{tpbdton102}  instead of a Gaussian bound as supposed in \cite{DOS} 
and the Avakumovic-Agmon-H\"ormander theorem for the spectral projection  of 
pseudo-differential operators on compact manifolds, then 
one argues as in Section 7.2 of  \cite{DOS}.
Even though the power of $\cos \theta$ in Theorem \ref{tpbdton102} is 
not optimal, it is then reduced by the interpolation argument as in the 
proofs of Theorems 3.1 or 3.2 in \cite{DOS}.
The advantage of  this method is that we obtain in addition that 
$f(A)$ is of weak type $(1,1)$ which is not stated in \cite{SeeSog}.

\appendix
\section{Compact manifolds} \label{Spbdtonapp}

Let $(M,g)$ be a Riemannian manifold (without boundary) of dimension $m$.
We always assume that a Riemannian manifold is $\sigma$-compact.
Then $M$ has a natural Radon measure denoted by $|\cdot|$. 
Let $p \in [1,\infty]$ and $k \in \Ni$.
Set 
\[
W^{k,p}_\loc(M)
= \{ u \in L_{p,\loc}(M) : u \circ \varphi^{-1} \in W^{k,p}_\loc(\varphi(V))
       \mbox{ for every chart } (V,\varphi) \} 
\;\;\; .  \]
If $u \in W^{1,p}_\loc(M)$ and $(V,\varphi)$ is a chart on $M$ then set 
$\frac{\partial}{\partial \varphi^i} u = (D_i (u \circ \varphi^{-1})) \circ \varphi \in L_{p,\loc}(V)$,
where $D_i$ denotes the partial derivative in $\Ri^m$.
Moreover, for all $u \in W^{1,p}_\loc(M)$, every chart $(V,\varphi)$ on $M$
and $i \in \{ 1,\ldots,m \} $ define $\nabla^i u, \nabla_i u \in L_{p,\loc}(V)$
by $\nabla^i u = \frac{\partial}{\partial \varphi^i} u$ and 
$\nabla_i u = \sum_{j=1}^m g_{ij} \, \nabla^j u$.
Note that $\nabla_i u$ and $\nabla^i u$ depend on the chart $(V,\varphi)$.
Let $k \in \Ni$ and $u \in W^{k,p}_\loc(M)$. 
Then there exists a unique element $|\nabla^k u| \in L_{p,\loc}(M)$
such that 
\[
|\nabla^k u| \Big|_V 
= \Big( \sum_{i_1,\ldots,i_k = 1}^m
   (\nabla_{i_1} \ldots \nabla_{i_k} u) \, 
   \overline{ (\nabla^{i_1} \ldots \nabla^{i_k} u) } \Big)^{1/2}
\]
for every chart $(V,\varphi)$ on $M$.
Set $|\nabla^0 u| = |u|$.
Similarly, if $u \in C^\infty(M \times M)$ and $k,\ell \in \Ni$, then 
there exists a unique element $|\nabla_{(1)}^k \, \nabla_{(2)}^\ell u| \in C(M \times M)$
such that 
\begin{eqnarray*}
|\nabla_{(1)}^k \, \nabla_{(2)}^\ell u| \Big|_V 
& = & \Big( \sum_{i_1,\ldots,i_k = 1}^m \sum_{j_1,\ldots,j_\ell = 1}^m
   (\nabla_{(1),i_1} \ldots \nabla_{(1),i_k} \nabla_{(2),j_1} \ldots \nabla_{(2),j_\ell} u) \, 
      \cdot \\*
& & \hspace*{50mm} {} \cdot
   \overline{ (\nabla_{(1)}^{i_1} \ldots \nabla_{(1)}^{i_k} 
               \nabla_{(2)}^{j_1} \ldots \nabla_{(2)}^{j_\ell} u) } \Big)^{1/2}
\end{eqnarray*}
for every chart $(V,\varphi)$ on $M$.
With obvious modifications one can also define 
$|\nabla_{(1)}^k \, \nabla_{(2)}^\ell u| \in C(M \times M)$ if $k = 0$ or $\ell = 0$.

Now also allow $k = 0$, so $k \in \Ni_0$.
Define the Banach space $W^{k,p}(M)$ by 
\[
W^{k,p}(M)
= \{ u \in W^{k,p}_\loc(M) : |\nabla^j u| \in L_p(M) \mbox{ for all } j \in \{ 0,\ldots,m \} \}
\]
with norm 
\[
\|u\|_{W^{k,p}(M)} = \Big( \sum_{j=0}^k \| \, |\nabla^j u| \, \|_p^2 \Big)^{1/2}
.  \]
If $u,v \in W^{1,2}(M)$ then there exists a unique
element $\nabla u \cdot \nabla v \in L_1(M)$ such that 
\[
(\nabla u \cdot \nabla v)|_V
= \sum_{i=1}^m (\nabla_i u) \, \overline{\nabla^i v}
\]
for every chart $(V,\varphi)$ on $M$.
Clearly if $(V,\varphi)$ is a chart on $M$ with $\overline V$ compact,
then for every multi-index $\gamma$ over $ \{ 1,\ldots,m \} $ 
there exists a $c > 0$ such that 
\begin{equation}
\Big\| \one_V \, \Big( \frac{\partial}{\partial \varphi} \Big)^\gamma u \Big\|_\infty
\leq c \, \|\nabla^{|\gamma|} u\|_\infty
\label{eSpbdtonapp;1}
\end{equation}
for all $u \in W^{|\gamma|,\infty}$.
Conversely, one has the following estimate on compact manifolds.

\begin{lemma} \label{lpbdtonA07}
Suppose $M$ is compact.
Let $N \in \Ni$ and for all $n \in \{ 1,\ldots,N \} $ let 
$(V_n,\varphi_n)$ be a chart on $M$ and $\chi_n \in C_c^\infty(V_n)$
such that $0 \leq \chi_n \leq \one$.
Suppose that $\sum_{n=1}^N \chi_n = \one$.
Let $k,\ell \in \Ni_0$.
Then there exists a $c > 0$ such that 
\[
|(\nabla_{(1)}^k \, \nabla_{(2)}^\ell u)(x,y)|
\leq c \sum_{n,m=1}^N
       \sum_{|\alpha| \leq k}
       \sum_{|\beta| \leq \ell}
   |\Big( 
\Big( \frac{\partial}{\partial \varphi_n} \Big)^\alpha_{(1)}
\Big( \frac{\partial}{\partial \varphi_m} \Big)^\beta_{(2)}
   ( (\chi_n \otimes \chi_m) u ) 
\Big) (x,y)|
\]
for all $u \in C^\infty(M \times M)$ and $x,y \in M$, 
where $( \frac{\partial}{\partial \varphi_n} )^\alpha_{(1)}$ acts on 
the first variable, we use multi-index notation, etc.
\end{lemma}

Define the sesquilinear form $a \colon W^{1,2}(M) \times W^{1,2}(M) \to \Ci$ by
$a(u,v) = \int \nabla u \cdot \nabla v$.
Then $a$ is closed and positive.
The {\bf Neumann Laplace--Beltrami operator} $\Delta$
on $M$ is the associated self-adjoint operator.
If $(V,\varphi)$ is a chart on $M$ then 
\[
\Delta \, u 
= \sum_{i,j=1}^d \frac{1}{\sqrt{g}} \, \frac{\partial}{\partial \varphi^i} \, 
           g^{ij} \, \sqrt{g} \, \frac{\partial}{\partial \varphi^j} \, u 
\]
for all $u \in C_c^\infty(V)$.
Since the form $a$ satisfies the Beurling--Deny criteria 
it follows that the semigroup $S$ generated by $\Delta$ extends to a 
continuous contraction semigroup $S^{(p)}$ on $L_p(M)$ for all $p \in [1,\infty)$.
We denote by $\Delta_p$ the generator of $S^{(p)}$.
If no confusion is possible, then we drop the suffix $p$ in $\Delta_p$.

\begin{prop} \label{ppbdtonA01} 
If $M$ is compact, $k \in \Ni$ and $p \in (1,\infty)$ then 
$W^{k,p}(M) = D((-\Delta_p)^{k/2})$.
Moreover, $C^\infty(M)$ is dense in $W^{k,p}(M)$.
\end{prop}
\proof\
See \cite{Hebey} Proposition~3.2.\hfill$\Box$

\vertspace

We need various Sobolev embeddings.

\begin{prop} \label{ppbdtonA02}
Suppose $M$ is compact.
Let $k,n \in \Ni_0$ and $p \in (2,\infty]$.
Suppose $\frac{1}{2} - \frac{1}{p} < \frac{k}{m}$.
Then $W^{k+n,2}(M) \subset W^{n,p}(M)$ and there
exists a $c > 0$ such that 
\[
\|u\|_{W^{n,p}(M)} \leq c \, \|u\|_{W^{k+n,2}(M)}^\alpha \, \|u\|_{L_2(M)}^{1 - \alpha}
\]
for all $u \in W^{k+n,2}(M)$, where $\alpha = \dfrac{n + m(\frac{1}{2} - \frac{1}{p})}{n+k}$.
\end{prop}
\proof\
These bounds are well known on $\Ri^m$ and then follow on a compact manifold
by localization.\hfill$\Box$

\begin{prop} \label{ppbdtonA02.2}
Suppose $M$ is compact.
Let $(V,\varphi)$ be a chart, $\chi \in C_c^\infty(V)$
and $\alpha$ a multi-index over $ \{ 1,\ldots,d-1 \} $.
Let $p \in (m,\infty)$ and $\ell \in \Ni$ be such that $\ell \geq |\alpha| + 1$.
Then there exists a $c > 0$ such that 
\[
\|\Big( \frac{\partial}{\partial \varphi} \Big)^\alpha (\chi \, u)\|_{L_\infty(M)}
\leq c \, \|\chi \, u\|_{W^{\ell,p}(M)}^\gamma \, 
            \|\chi \, u\|_{L_\infty(M)}^{1 - \gamma}
\]
for all $u \in C^\infty(M)$, where 
\[
\gamma = \frac{|\alpha|}{\ell - \frac{m}{p}}
.  \]
\end{prop}
\proof\
By the Sobolev embedding theorem and interpolation there exists a $c' > 0$ 
such that 
\[
\|v\|_{W^{|\alpha|,\infty}(\Ri^m)}
\leq c' \, \|v\|_{W^{\ell,\infty}(\Ri^m)}^\gamma \, \|v\|_{L_\infty(\Ri^m)}^{1-\gamma}
\]
for all $v \in W^{\ell,p}(\Ri^m)$.
Using the chart $(V,\varphi)$ and localizing with $\chi$ gives the 
proposition.\hfill$\Box$

\begin{lemma} \label{lpbdtonA02.5} 
Let $(V_1,\varphi)$ and $(V_2,\psi)$ be charts on $M$, 
let $\chi_1,\chi_2 \in C_c^\infty(M)$ and suppose that 
$\supp \chi_1 \subset V_1$ and $\supp \chi_2 \subset V_2$.
Let $k \in \Ni_0$ and $T \in \OPS^k(M)$.
Let $p \in (1,\infty)$.
Then for every multi-index $\alpha$ over $ \{ 1,\ldots,m \} $ with 
$|\alpha| \leq k$ there exists a bounded operator $T_\alpha$ on 
$L_p(M)$ such that 
\[
M_{\chi_1} \, T \, M_{\chi_2}
= \sum_{|\alpha| \leq k} 
    M_{\chi_1} \, T_\alpha \, 
       \Big( \frac{\partial}{\partial \psi} \Big)^\alpha \, M_{\chi_2}
.  \]
\end{lemma}
\proof\
There exists a $\widetilde T \in \OPS^k(\Ri^m)$ such that 
\[
\widetilde T w 
= \Bigg( \chi_1 \, T \bigg( \Big( w \cdot (\chi_2 \circ \psi^{-1}) \Big) 
           \circ \psi \bigg) \Bigg) \circ \varphi^{-1}
\]
for all $w \in \cs(\Ri^m)$.
By the proof of Proposition~VI.5 in \cite{Ste3} for all
multi-indices $\alpha$ with $|\alpha| \leq k$ there exists a 
pseudo-differential operator $\widetilde T_\alpha$ of order $0$ such that 
$\widetilde T = \sum_{|\alpha| \leq k} \widetilde T^\alpha \in \partial^\alpha$.
Each $\widetilde T_\alpha$ is bounded on $L_p(\Ri^d)$ by \cite{Ste3} Proposition~VI.4.
Then the lemma follows by a coordinate transformation.\hfill$\Box$

\begin{lemma} \label{lpbdtonA03} 
Suppose $M$ is compact.
Let $k \in \Ni_0$ and $T \in \OPS^k(M)$.
Let $p \in (1,\infty)$.
Then there exists a $c > 0$ such that 
$\|T u\|_p \leq c \, \|u\|_{W^{k,p}(M)}$
for all $u \in C^\infty(M)$.
\end{lemma}
\proof\
This follows with a partition of the unity from Lemma~\ref{lpbdtonA02.5}.\hfill$\Box$

\begin{lemma} \label{lpbdtonA03.5} 
Let $(V,\varphi)$ be a chart on $M$ and $\chi \in C_c^\infty(V)$.
Let $k \in \Ni_0$ and $T \in \OPS^k(M)$.
Let $p \in (1,\infty)$.
Then for every multi-index $\alpha$ over $ \{ 1,\ldots,m \} $ with 
$|\alpha| \leq k$ there exists a bounded operator $T_\alpha$ on 
$L_p(M)$ such that 
\[
M_\chi \, T
= \sum_{|\alpha| \leq k} 
    M_\chi \, 
       \Big( \frac{\partial}{\partial \varphi} \Big)^\alpha \, T_\alpha
.  \]
\end{lemma}
\proof\
This follows from Lemma~\ref{lpbdtonA02.5}, duality and a partition of the 
unity.\hfill$\Box$

\vertspace

For the remaining part of this section suppose that the manifold $M$ is connected.
Then the Riemannian manifold has a natural distance, denoted by~$d_M$.
Note that 
\begin{equation}
d_M(x,y) 
= \sup \{ |g(x) - g(y)| : g \in C^\infty(M) \mbox{ and } \|\nabla g\|_\infty \leq 1 \} 
\label{elpbdtonA05;1}
\end{equation}
for all $x,y \in M$.
See, for example \cite{ABE} Proposition~2.2.
We need some equivalence of the distance on $M$.
Since $M$ is compact, one can locally regularize using a finite number of 
charts. 
Therefore (\ref{elpbdtonA05;1}) implies the next lemma.

\begin{lemma} \label{lpbdtonA05}
For all $N \in \Ni$ there exists a $c > 0$ such that 
\begin{eqnarray*}
\frac{1}{c} \, d_M(x,y)
& \leq & \sup \{ g(x) - g(y) : g \in C^\infty(M) \mbox{ and }
      \|\nabla^i g\|_\infty \leq 1 \mbox{ for all } i \in \{ 1,\ldots,N \} \}  \\
& \leq & c \, d_M(x,y)
\end{eqnarray*}
for all $x,y \in M$.
\end{lemma}

Moreover, for embedded manifolds the distance $d_M$ is comparable 
with the Euclidean distance.
This is a consequence of \cite{Hel2} Proposition~9.10.

\begin{lemma} \label{lpbdtonA06}
Suppose $k \in \Ni$ and $M$ is embedded in $\Ri^k$.
Then there exists a $c > 0$ such that 
\[
\frac{1}{c} \, d_M(x,y)
\leq |x - y| 
\leq c \, d_M(x,y)
\]
for all $x,y \in M$.
\end{lemma}

Finally we introduce H\"older spaces.
If $\nu \in (0,1)$ then we denote by $C^\nu(M)$ the space of all H\"older continuous
functions of order $\nu$ with respect to the distance~$d_M$, with seminorm
\[
|||u|||_{C^\nu(M)} 
= \sup_{x \neq y} \frac{|u(x) - u(y)|}{d_M(x,y)^\nu}
.  \]
The norm on $C^\nu(M)$ is given by 
$\|u\|_{C^\nu(M)} = \|u\|_\infty + |||u|||_{C^\nu(M)}$.
With this norm the space $C^\nu(M)$ is a Banach space.
Moreover, one has the following Sobolev embedding.

\begin{prop} \label{ppbdtonA04}
Suppose $M$ is compact and $p \in (m,\infty)$.
Set $\nu = 1 - \frac{m}{p}$.
Then $W^{1,p}(M) \subset C^\nu(M)$.
In particular, there exists a $c > 0$ such that 
\[
\|u\|_{C^\nu(M)} \leq c \, \|u\|_{W^{1,p}(M)}
\]
for all $u \in W^{1,p}(M)$.
\end{prop}
\proof\
See \cite{Hebey} Theorem~3.5.\hfill$\Box$

\subsection*{Acknowledgements} 
The authors wish to thank Wolfgang Arendt, Gerd Grubb and J\"org Seiler 
for useful discussions.
The main part of this work was carried out whilst  the first named
author was visiting the  Institut  of Mathematics at the University of  Bordeaux 1  
in September 2012.  
He wishes to thank  the University of Bordeaux 1 for financial support. 
The research of A.F.M. ter  Elst  is partly supported by the 
Marsden Fund Council from Government funding, 
administered by the Royal Society of New Zealand. 
The research of E.M.  Ouhabaz  is partly supported by the ANR 
project `Harmonic Analysis at its Boundaries',  ANR-12-BS01-0013-02.

\end{document}